\documentstyle[12pt,leqno,amssymb,amscd,graphics,hyperref]{amsart}  
\hypersetup{hidelinks=true}

\newtheorem{thm}{Theorem}[section]
\newtheorem{lemma}[thm]{Lemma}
\newtheorem{prop}[thm]{Proposition}
\newtheorem{cor}[thm]{Corollary}
\newtheorem{clcor}[thm]{Classical Corollary}
\theoremstyle{definition}
\newtheorem{dfn}[thm]{Definition}
\theoremstyle{remark}
\newtheorem{remark}[thm]{Remark}
\newtheorem{example}[thm]{Example}

\begin{document}

\newcommand{\ct}{\cite}
\newcommand{\pr}{\protect\ref}
\newcommand{\su}{\subseteq}
\newcommand{\pa}{{\partial}}
\newcommand{\e}{\epsilon}
\newcommand{\es}{{\varnothing}}

\newcommand{\DD}{{\lambda}}

\newcommand{\QQ}{{\mathbb Q}}
\newcommand{\RR}{{\mathbb R}}
\newcommand{\rn}{{{\mathbb R}^n}}
\newcommand{\ZZ}{{\mathbb Z}}
\newcommand{\N}{{\mathbb N}}

\newcommand{\h}{\hat}
\newcommand{\dd}{\delta}
\newcommand{\x}{\Delta}
\newcommand{\xxk}{{\x^k_{v_1,\dots,v_k}}}
\newcommand{\xt}{{\x^2_{v,w}}}
\newcommand{\sm}{\s  \! M}
\newcommand{\sx}{\s  \! X}

\newcommand{\1}{^{1 \over \DD}_{\phantom i}}

\newcommand{\s}{{}^*}

\newcommand{\f}{{\varphi}}

\newcommand{\ha}{{\mathfrak{h}}}
\newcommand{\cp}{{\mathcal{P}}}
\newcommand{\cu}{{\mathcal{U}}}
\newcommand{\cd}{{\mathcal{D}}}
\newcommand{\cs}{{\mathcal{S}}}
\newcommand{\cc}{{\mathcal{C}}}
\newcommand{\lf}{{\lfloor}}
\newcommand{\rf}{{\rfloor}}
\newcommand{\vv}{\overrightarrow}
\newcommand{\tf}{{\widetilde{F}}}
\newcommand{\tg}{{\widetilde{G}}}
\newcommand{\der}{{\partial}}

\newcommand{\tth}{{   {}^{\mathfrak{h}}     }}

\newcommand{\tm}{{{}^{\mathfrak{h}} \! {M}}}

\newcommand{\ty}{{{}^{\mathfrak{h}}  {Y}}}
\newcommand{\tu}{{{}^{\mathfrak{h}}  {U}}}
\newcommand{\tr}{{{}^{\mathfrak{h}}  {\RR}}}
\newcommand{\ii}{{(0,1)}}

%\subject{primary}{msc2000}{26E35, 58A05}
%\subject{secondary}{msc2000}{}

\newcounter{numb}

\title{Differential Geometry via Infinitesimal Displacements}

\author{Tahl Nowik}
%\givenname{Tahl}
%\surname{Nowik}
\address{Department of Mathematics, Bar-Ilan University, 
Ramat-Gan 5290002, Israel}
\email{tahl@@math.biu.ac.il}
\urladdr{www.math.biu.ac.il/~tahl}

\author{Mikhail G. Katz}
%\givenname{Mikhail G.}
%\surname{Katz}
\address{Department of Mathematics, Bar-Ilan University, 
Ramat-Gan 5290002, Israel}
\email{katzmik@@math.biu.ac.il}
\urladdr{www.math.biu.ac.il/~katzmik}

\begin{abstract}
We present a new formulation of some basic differential geometric notions on a smooth manifold $M$, in the setting of nonstandard analysis.
In place of classical vector fields, for which one needs to construct the tangent bundle of $M$, we define  a \emph{prevector field}, which is an internal map from $\sm$ to itself,
implementing the intuitive notion of vectors as infinitesimal displacements. 
We introduce regularity conditions for prevector fields, defined by
finite differences, thus purely combinatorial conditions involving no
analysis.  These conditions replace the more elaborate analytic
regularity conditions appearing in previous similar approaches,
e.g. by Stroyan and Luxemburg or Lutz and Goze.
We define the flow of a prevector field by hyperfinite iteration of the given prevector field, in the spirit of Euler's method. We define the Lie bracket of two prevector fields by appropriate iteration of their commutator. 
We study the properties of flows and Lie brackets, particularly in relation with our proposed regularity conditions. We present several simple applications
to the classical setting, such as bounds related to the flow of vector fields, analysis of small oscillations of a pendulum, and an instance  of Frobenius' Theorem regarding the complete integrability of independent vector fields.
\end{abstract}

\maketitle

\section{Introduction}\label{intro}

We develop foundations for differential geometry  on  smooth manifolds,
based on infinitesimals, where vectors and vector fields are represented 
by infinitesimal displacements in the manifold itself, as they were thought of historically.
Such an approach was previously introduced e.g. by Stroyan and Luxemburg in \ct{stl}, and by 
Lutz and Goze in \ct{lg}. For such an approach to work, one needs to assume some regularity
condition on the infinitesimal displacement maps used to represent vector fields, in place of the smoothness properties appearing in the classical setting.
The various regularity conditions chosen in existing sources seem non basic and overly tied up with 
classical analytic notions. In the present work we  introduce natural and easily verifiable regularity conditions, defined by finite differences. 
We show that these weak regularity conditions are sufficient for defining notions such as the flow and Lie bracket of vector fields.

In more detail, we would like to study vectors and vector fields in a smooth manifold $M$ while bypassing its tangent bundle. We would like to think
of a vector based at a point $a$ in $M$ as a little arrow in $M$ itself, whose tail is $a$ and whose head is a nearby point $x$ in $M$. The notions ``little'' and ``nearby'' can be formalized in the hyperreal framework, i.e. nonstandard analysis, where \emph{infinitesimal} quantities are available. 
We aim to demonstrate that various differential  geometric concepts and various proofs become simpler and more transparent in this framework. 
For good introductions to nonstandard analysis see e.g. Albeverio et al \ct{a}, Goldblatt \ct{g}, Gordon, Kusraev and Kutateladze \ct{GKK}, Keisler \ct{k}, Loeb and Wolff \ct{lw}, V\"{a}th \ct{v}. For an advanced study of axiomatic treatments see Kanovei and Reeken \ct{kr}. For a historical perspective see Bascelli et al \ct{bb}.
For additional previous applications of nonstandard analysis to differential geometry see Almeida, Neves and Stroyan \ct{ans} and references therein.
For application of nonstandard analysis to the solution of Hilbert's fifth problem see Tao \ct{tao}  and references therein. Nonstandard analysis was initiated by Robinson~\ct{rob}.

Infinitesimal quantities may themselves be infinitely large or infinitely small compared to one another, so the key to our application of infinitesimal quantities in differential geometry is to fix a positive infinitesimal hyperreal number $\DD$ once and for all, which will fix the scale of our constructions.
We then define a \emph{prevector} based at a nearstandard point $a$ of $\sm$ to be a pair of points $(a,x)$ in $\sm$ for which
the distance between $a$ and $x$ is not infinitely large compared to $\DD$. Two prevectors $(a,x)$, $(a,y)$ based at  $a$ are termed equivalent if the distance between $x$ and $y$ is
infinitely small compared to $\DD$. Note that the notion of the distance in $\sm$ being infinitely large or infinitely small compared to $\DD$ does not require a metric on $M$; it is intrinsic to the differentiable structure of $M$.
We next define a prevector \emph{field} to be an internal map $F:\sm \to \sm$ such that for every nearstandard point $a$ in $\sm$, the pair $(a,F(a))$ is  a prevector at $a$. The requirement that the map
 $F$ be internal is crucial e.g. for hyperfinite iteration, internal induction, and the internal definition principle.    
Two prevector fields $F,G$ are equivalent if for every nearstandard $a$ in $\sm$, the pairs $(a,F(a))$ and $(a,G(a))$ are equivalent prevectors.

As already mentioned,
in place of the hierarchy $C^k$ of smoothness appearing in the classical setting of vector fields, we introduce a hierarchy $D^k$
of weaker regularity conditions for prevector fields, defined by finite differences. 
We show that these weaker regularity conditions are sufficient for defining notions such as the flow of a prevector field and 
the Lie bracket of two prevector fields.
We use our flow to show that a canonical representative can be chosen from every equivalence class of local prevector fields,
and that every classical vector field on $M$ can be realized by a prevector field on $\sm$ whose values on the nearstandard part of $\sm$ are canonically prescribed.   
For this last statement, in case the manifold is non compact, we will need to  
assume that our nonstandard extension is countably saturated.  (The notion of countable saturation will be explained in Section~\pr{rl}.)

The framework we propose suggests various possibilities for further
investigation. 
For example, one could seek to formulate notions
corresponding to the Poincar\'{e}-Hopf Theorem regarding indices of zeros of vector fields, in terms of infinitesimal displacements.
Another example is proving Frobenius' Theorem in the spirit of our proof of Theorem~\pr{comm} and 
Classical~Corollary~\pr{cc4}, characterizing when $k$ vector fields are the first $k$ coordinate vector fields for some choice of coordinates. This can be thought of as an 
instance of Frobenius' Theorem, with stronger assumption and stronger conclusion. 

A very different  alternative approach to the foundations of differential geometry is that of Synthetic Differential Geometry, introduced by Lawvere and others, see e.g. Kock \ct{ko}. It relies on category-theoretic concepts and intuitionistic logic, which are not needed for our approach.
To the extent that our hierarchy $D^k$ of regularity
classes is formulated in terms of finite differences and thus avoids
classical \emph{analytic} notions, our approach can also be characterized as
\emph{synthetic} differential geometry.

The structure of the paper is as follows.
In Section~\pr{tv} 
we define prevectors, and explain their relation to classical tangent vectors. We define such notions as 
the action of a prevector on a smooth function, and the differential of a smooth map from one smooth manifold to another.
In Section~\pr{vf} 
we define local and global prevector fields. We define the $D^k$ regularity property of prevector fields.
The property $D^k$ is defined via coordinates by a finite difference condition.  
We show how a local classical vector field induces a local prevector field,
and show that if the classical vector field is $C^k$ then the induced local prevector field is $D^k$ (Proposition~\pr{6}).
We then show the following,
which though rather simple for $D^1$ turns out somewhat involved for $D^2$.
\begin{thm}
The definition of $D^1$ and $D^2$ prevector fields is independent of the choice of coordinates (Propositions \pr{5}, \pr{4}).
\end{thm}
We further show that $D^2$ implies $D^1$ (Proposition~\pr{3}).  
In Section~\pr{global} we show that a global $D^1$ prevector field is bijective on the nearstandard part of $M$ (Theorem~\pr{bij}).   

In Section~\pr{f} 
we define the flow of a prevector field by hyperfinite iteration of the given prevector field. 
It is a generalization of the Euler approximation for the flow appearing e.g. in Keisler \ct[p.~162]{k}, as well as  Stroyan and Luxemburg \ct[p.~128]{stl}  and Lutz and Goze \ct[p.~115]{lg}. Using straightforward internal induction we prove the following.
\begin{thm} 
The flow of a $D^1$ prevector field remains in a bounded region for some appreciable interval of time.
The growth of the distance between two points moving under the flow  
is bounded above and below after time $t$ by factors of the form $e^{\pm K t}$ (Theorem~\pr{flow1}).  The difference between the flows of different prevector fields is bounded above by a function of the form $\beta t e^{K t}$ (Theorem~\pr{flow2}). 
\end{thm}

We note that this implies corresponding bounds for the flow of classical vector fields
(Classical~Corollary~\pr{cc1}).
 We use the flow of prevector fields to show the following.
\begin{thm}
A canonical representative can be chosen from each equivalence class of local prevector fields which contains a $D^1$ (resp. $D^2$)  prevector field, and this 
representative is itself $D^1$ (resp. $D^2$).
\end{thm}
This is done roughly as follows.  Given a $D^1$ (resp. $D^2$) local prevector field $F$, its flow in $\sm$ induces a standard local flow in $M$, which is then extended back to $\sm$ and evaluated at time $t=\DD$, producing a new prevector field $\tf$. 
Different representatives $F$ of the given equivalence class induce the same standard flow (Theorem~\pr{hft}) and so $\tf$ is indeed 
canonically chosen.  We then need to show that $\tf$ 
is in fact equivalent to the original $F$ one started with (Theorem~\pr{tld}), and that if $F$ is
$D^1$ or $D^2$ then the same holds for $\tf$ (Propositions \pr{111}, \pr{112}).    
As example of an application of our results on flows we 
analyze oscillations of a pendulum with infinitesimal amplitude (Section~\pr{secpend}). 

In Section~\pr{rl} 
we show the following.
\begin{thm}
Every \emph{global} classical $C^1$ (resp. $C^2$) vector field can be realized by a \emph{global} $D^1$ (resp. $D^2$) prevector field, whose values on the nearstandard part of $\sm$ are canonically prescribed (Theorem~\pr{gl}). 
\end{thm}
This involves techniques similar to those mentioned above in relation to the construction of $\tf$, and additionally, for a vector field which does not have compact support, our assumption of countable saturation is used. 
 
 In Section~\pr{lbr} 
we define the Lie bracket of two prevector fields, and relate it to the classical 
Lie bracket of classical vector fields (Theorem~\pr{rb}).  We show the following.
\begin{thm}
The Lie bracket of two $D^1$ prevector fields is itself a prevector field (Theorem~\pr{7}). The Lie bracket of two $D^2$ prevector fields is $D^1$ (Theorem~\pr{9}).
The Lie bracket is well defined on equivalence classes of $D^2$ prevector fields (Theorem~\pr{ld2}),
 and this is not the case for prevector fields that are merely $D^1$ (Example~\pr{e4}).
\end{thm}
We show that the Lie bracket of two $D^2$ prevector fields is equivalent to the identity prevector field if and only if their local standard flows commute (Theorem~\pr{comm}). We note that this implies the classical 
result that the flows of two vector fields commute if and only if their Lie bracket vanishes (Classical~Corollary~\pr{cc4}).

We would like to thank Thomas McGaffey for guiding us to the existing literature on nonstandard analysis approaches to differential geometry.

\section{Prevectors}\label{tv}

For our analysis we will need to compare different infinitesimal quantities. 
It is helpful to introduce the following relations.

\begin{dfn}\label{dd1}
For $r,s \in \s\RR$, we will write $r \prec s$ if $r=as$ for finite $a$, and will write  $r \prec\prec s$ if  $r=as$ for infinitesimal $a$.
\end{dfn}

Thus $r\prec 1$ means that $r$ is finite, and $r \prec\prec 1$ means that $r$ is infinitesimal.
Given a finite dimensional vector space $V$ over $\RR$,  and given $v \in \s V$, and $s\in \s\RR$, we will write $v\prec s$ if 
$\s\|v\| \prec s$ for some norm $\|\cdot\|$ on $V$. 
We will generally omit the $*$ from function symbols and so will simply write $\|v\| \prec s$.  
This condition is independent of the choice of norm since all norms on $V$ are equivalent. Similarly we will write $v\prec\prec s$ if $\|v\| \prec\prec s$.
We will also write $v \approx w$ when $v-w \prec\prec 1$.
If one chooses a basis for $V$ thus identifying it with $\RR^n$, then $x=(x_1,\dots,x_n)\in\s\RR^n$ satisfies $x\prec s$ or $x\prec\prec s$     if and only if  each $x_i$ satisfies this, since this is clear in, say, the Euclidean norm on $\RR^n$. 
Given $0 < s\in\s\RR$ 
let $V^s_F=\{v\in \s V:v \prec s\}$ and  $V^s_I=\{v\in \s V:v \prec\prec s\}$.  Then $V^s_I \su V^s_F \su \s V$ are linear subspaces over $\RR$, and  $V^s_F / V^s_I \cong V$, since
it is well known that $V^1_F /V^1_I \cong V$, and 
multiplication by $s$ maps $V^1_F$ onto $V^s_F$ 
and $V^1_I$ onto $V^s_I$.

Our object of interest is a smooth manifold $M$. For $p \in M$, the halo of $p$, which we denote  by $\ha(p)$, is the set of all points $x$ in the nonstandard extension $\sm$ of $M$,
for which  there is a coordinate neighborhood $U$ of $p$ such that $x\in\s U$ 
and $x \approx p$ in the given coordinates. In fact, the definition of $\ha(p)$ does not require coordinates, but rather depends only on the \emph{topology} of $M$.
\footnote{For a topological space $X$ and $p \in X$ let $N_p$ be the set of all open neighborhoods of $p$ in $X$. Then
the halo (or monad) of $p$ is defined as  $\ha (p)= \bigcap_{U \in N_p} \s U$.}
The points of $M$ are called \emph{standard}, 
and a point which is in $\ha(p)$ for some $p\in M$ 
is called \emph{nearstandard}. If $a$ is nearstandard then  the standard part (or shadow) of $a$, denoted $st(a)$, is the unique $p \in M$ such that $a \in \ha(p)$.

\begin{dfn}\label{dd2}
For $A\su M$, the halo of $A$ is  $\tth \! A = \bigcup_{a\in A} \ha(a)$. 
\end{dfn}

In particular, $\tm$ is the set of all nearstandard points in $\sm$.
If  $A\su M$ is open in $M$ then $\tth\! A \su \s\! A$,
and if it is compact then $ \s\! A \su   \tth\! A $. 
In particular, if $M$ is compact then $\tm = \sm$.    
If $M$ is noncompact then $\tm$ is an external set. 
\footnote{The converse of the statements in this paragraph also hold if we assume that our nonstandard extension satisfies countable saturation, and using the fact that $M$ has a countable basis. See e.g. Albeverio et al \ct[Section 2.1]{a}}

Much of our analysis will be local, so
given an open $W \su \RR^n$ and a smooth function $f:W \to \RR$ we 
note some properties of the extension $\s\! f : \s W \to \s \RR$, obtained by transfer. 
When there is no risk of confusion, we will omit the $*$ from the function symbol $\s\! f$ and simply write $f$ for both the original function and its extension.

\begin{lemma}\label{ff}
For open $W\su\rn$,
let $f:W \to\RR$ be continuous. Then $f(a)$ is finite for every $a \in \tth W$. 
\end{lemma}

\begin{pf}
Given $a \in \tth W$,
let $U$ be a neighborhood of $st(a)$ such that 
$\overline{U} \su W$ and $\overline{U}$  is compact. So there is $C \in \RR$ such that $|f(x)| \leq C$ for
all $x \in U$. 
By transfer $|f(x)| \leq C$ for all $x \in \s U$, in particular $|f(a)| \leq C$, so $f(a)$ is finite.
(As for functions, we omit the $*$ from relation symbols, writing $\leq$ in place of $\s\!\!\leq$.)
\end{pf}

Given an open $U \su \RR^n$ and a smooth function $f:U \to \RR$,
the partial derivatives of $\s\! f$ are by definition the functions $\s\big( {\pa f \over \pa x_i}\big)$, i.e.
the extensions of the partial derivatives of $f$.
So, one has a row vector $D_a$ of partial derivatives at any point $a\in\s U$ if $f$ is $C^1$, and similarly a Hessian matrix $H_a$ of the second partial derivatives at every $a \in\s U$,
in case $f$ is $C^2$.
By Lemma~\pr{ff} $D_a$  and $H_a$ are finite throughout $\tth U$.

We state the following properties of $D_a$ and $H_a$ as three remarks for future reference.

\begin{remark}\label{dh1}
Let $a,b \in \tu$ with $a\approx b$, then the interval between $a$ and $b$ is included in $\tu \su \s U$.
By transfer of the mean value theorem, if $f$ is $C^1$ then $f(b) - f(a) = D_x(b-a)$ for some $x$ in the interval between $a$ and $b$.
Since the partial derivatives are continuous, we have by the characterization of continuity via infinitesimals 
\footnote{$g$ is continuous at $a$ if and only if  $\s\! g(\ha(a)) \su \ha(g(a))$.} 
that
$D_x - D_y \prec\prec 1$ for any $y \approx a$ (e.g. $y=st(a)$), 
so writing $$f(b)-f(a) = D_y(b-a) + (D_x-D_y)(b-a),$$ we see that $f(b)-f(a) - D_y(b-a) \prec\prec \|b-a\|$. 
Furthermore, if the first partial derivatives are Lipschitz in some neighborhood  (e.g. if $f$ is $C^2$), 
and we are given a constant $\beta\prec\prec 1$ such that
$\|x - y \| \prec \beta$ for all $x$ in the interval between $a$ and $b$,
then we have the stronger condition $f(b)-f(a) - D_y(b-a) \prec \beta\|b-a\|$. 
\end{remark}

\begin{remark}\label{dh2}
If $f$ is $C^2$ then by transfer of the Taylor approximation theorem we have
$f(b)-f(a) = D_a(b-a) +{1\over 2} (b-a)^tH_x(b-a)$ for some $x$ in the interval between $a$ and $b$, and remarks similar to those we have made regarding $D_x - D_y$ apply to $H_x - H_y$.
\end{remark}

\begin{remark}\label{dh3}
If $\f = ( \f^i ) : U \to \RR^n$ then the 
$n$ rows $D^i_a$ corresponding to $\f^i$ form the Jacobian matrix $J_a$ of 
$\f$ at $a$. By applying the above considerations to each $\f^i$ we obtain that 
$\f(b)-\f(a) - J_y(b-a) \prec\prec \|b-a\|$, or if all partial derivatives are Lipschitz
in some neighborhood 
 (e.g. if $\f$ is $C^2$) then $\f(b)-\f(a)  - J_y(b-a) \prec \beta \|b-a\|$ with $\beta$ as above.
\end{remark}

Now choose a positive infinitesimal $\DD \in \s{\RR}$, and \emph{fix it once and for all}.

\begin{dfn}\label{p}
Given  $a \in \tm$,   a prevector based at $a$ is a pair $(a,x)$, $x \in \tm$, such that  for every smooth $f:M \to \RR$,
$f(x) - f(a) \prec \DD$. 
Equivalently, given coordinates in a neighborhood $W$ of $st(a)$ in $M$, whose image is  $U \su \RR^n$,
 and $\hat{a},\hat{x} \in \s U$ are the coordinates for $a$ and $x$, then
$(a,x)$ is a prevector based at $a$ if $\hat{x} - \hat{a} \prec \DD$, where the difference   $\hat{x} - \hat{a}$ is defined in $\s\RR^n \supseteq \s U$.
\end{dfn}

We show that the two definitions are indeed equivalent. Assume the first definition, and let $x_1,\dots,x_n$ be the chosen coordinate functions. Since each $x_i$ is smooth, 
we get  $ x_i(x)-  x_i(a) \prec \DD$ for each $i$, i.e.  $\hat{x} - \hat{a} \prec \DD$. 
Conversely, assume the second definition holds, and let $f:M \to \RR$ be a smooth function. Then  
$ f(x) -  f(a) = D_c(\hat{x}-\hat{a})$ for some $c$ in the interval between $\hat{a}$ and $\hat{x}$, and so 
$f(x)-f(a) \prec \| \hat{x} - \hat{a}\| \prec \DD$.
(The components of $D_c$ are finite by Lemma~\pr{ff}.)

We denote by $P_a=P_a(M)$ the set of prevectors based at $a$.

\begin{dfn}\label{dd3}
We define an equivalence relation $\equiv$ on $P_a$ as follows: $(a,x)\equiv (a,y)$ if $ f(y) -  f(x) \prec\prec\DD$ for 
 every smooth $f:M \to \RR$, or equivalently, if in coordinates as above, $\h{y}-\h{x} \prec\prec\DD$. 
\end{dfn}

The equivalence of the two definitions follows by the same argument as above.
Since the relation $(a,x)\equiv (a,y)$ depends only on $x,y$, we will also simply write $x \equiv y$.
We denote the set of equivalence classes $P_a/\!\!\equiv$ by $T_a=T_a(M)$.
In the spirit of physics notation, the equivalence class of 
$(a,x) \in P_a$ will be denoted $\vv{ax}$.

Given $a \in \tm$  let $W$ be a coordinate neighborhood of $st(a)$ in $M$ with image $U \su\RR^n$.
We can identify $P_a$ with  $(\RR^n)^\DD_F$ via $(a,x) \mapsto \h{x}-\h{a}$. This induces an identification of $T_a = P_a/\!\!\equiv$ 
with $\RR^n = (\RR^n)^\DD_F/(\RR^n)^\DD_I$.
Under this identification $T_a$ inherits the structure of a vector space over $\RR$.
If we choose different coordinates in a neighborhood of $st(a)$ with image $U' \su \RR^n$, 
then if $\varphi:U \to U'$ is the change of coordinates, 
then by Remark~\pr{dh3}  $\varphi(\h{x}) - \varphi(\h{a}) - J_{st(a)}(\h{x} - \h{a}) \prec\prec \|\h{x} - \h{a}\|\prec \DD$.
This means that the map $\RR^n \to \RR^n$ induced by the two identifications of $T_a$ with $\RR^n$ provided by the two coordinate maps, is given by multiplication by the matrix $J_{st(a)}$, and so is linear.
Thus the vector space  structure induced on $T_a$ via coordinates is independent of the choice of coordinates, and so we have a well defined 
vector space structure on $T_a$ over $\RR$.
Note that the object $T_a$ is a mixture of standard and nonstandard notions. It is a vector space over
$\RR$ rather than $\s\RR$, but defined at every $a\in \tm$.

If $a,b \in \tm$ and $a\approx b$, then given coordinates in a neighborhood of $st(a)$, the identifications of $T_a$ and $T_b$ with $\RR^n$ induced by these coordinates induces an identification between $T_a$ and $T_b$. Given a different choice of coordinates, the matrix $J_{st(a)}$ used in the previous 
paragraph is the same matrix for $a$ and $b$,  and so the identification of $T_a$ with $T_b$ is well defined, independent of a choice of coordinates. Thus  when $a \approx b \in \tm$ we may unambiguously add a vector 
$\vv{a \ x} \in T_a$ with a vector $\vv{b \ y} \in T_b$.

\begin{dfn}\label{dd4}
A prevector $(a,x) \in P_a$ acts on a smooth function $f:M \to \RR$ as follows:   
$$(a,x)f = {1\over\DD}( f(x)- f(a)),$$ which is finite by definition of prevector.
\end{dfn}

This induces a differentiation of $f:M \to \RR$ by a vector $\vv{a \ x} \in T_a$ as follows:   
$\vv{a \ x} \ f = st((a,x)f)$. The action $\vv{a \ x} \ f$ 
is well defined by definition of the equivalence relation $\equiv$.
Note our mixture again, $\vv{a \ x}$ is a nonstandard object based at the nonstandard point $a$, but it assigns a standard real number to the standard function $f$.
The action $(a,x)f$  satisfies the Leibniz rule up to infinitesimals, indeed:
\begin{align*}
{1\over\DD}\Bigl( f(x)g(x)- f(a)g(a) \Bigr) & = {1\over\DD}\Bigl( f(x)g(x)-f(x)g(a)+ f(x)g(a) -f(a)g(a) \Bigr)  \\
& =f(x){1\over\DD}\Bigl(g(x)-g(a)\Bigr)+{1\over\DD} \Bigl(f(x) -f(a)\Bigr)g(a)  \\  & \approx
f(a){1\over\DD}\Bigl(g(x)-g(a)\Bigr)+{1\over\DD} \Bigl(f(x) -f(a)\Bigr)g(a)
\end{align*} 
where the final $\approx$ is by continuity of $f$.
For the action $\vv{a \ x}f$ this implies the following,  where the second equality is by continuity of $f$ and $g$. 

\begin{prop}\label{pvec}
Letting $a_0=st(a)$ we have
$$\vv{a \ x}(fg) = st(f(a)) \cdot \vv{a \ x} g + \vv{a \ x}f \cdot st(g(a))= f(a_0) \cdot \vv{a \ x}g + \vv{a \ x}f \cdot g(a_0).$$ 
\end{prop}

\begin{dfn}\label{dd5}
If $h:M \to N$ is a smooth map between smooth manifolds, then for $a \in \tm$ we define the differential of $h$, 
$dh_a:P_a(M) \to P_{h(a)}(N)$   by setting $$dh_a((a,x)) = ( h(a)  , h(x)).$$ 
\end{dfn}

This induces a map 
$dh_a:T_a(M) \to T_{h(a)}(N)$ given by $dh_a(\vv{a \ x}) = \vv{ h(a)  \  h(x)}$. 
The relation here between $h$ and $dh$ seems more transparent than in the corresponding classical definition. Furthermore, the ``chain rule'', i.e. the fact that $d({g \circ h})_a = dg_{h(a)} \circ dh_a$, becomes immediate, for both $P_a$ and $T_a$. Namely,  
$dg_{h(a)} \circ dh_a((a,  x)) =dg_{h(a)} (( h(a),  h(x))) =
 ( g\circ h(a)  , g \circ h(x)) = d({g \circ h})_a((a , x))$, and
similarly for $\vv{a \ x} \in T_a$.

\begin{remark}\label{67}
For standard $a \in M$, $T_a$ is naturally identified with the classical tangent space of $M$ at $a$ as follows. Since we have done everything also in terms of coordinates, it is enough to see this 
for open $U \su \RR^n$, where the tangent space at any point $a$ is $\RR^n$ itself. A vector $v \in \RR^n$  is then identified with 
 $\vv{a \ (a+\DD\cdot v)}$. Under this identification, our definitions of $\vv{a \ x}f$ and  $dh(\vv{a \ x})$ coincide with the classical ones. 
\end{remark}

\section{Prevector fields}\label{vf}

For a smooth manifold $M$,
recall that $\tm$ denotes the set of all nearstandard points in $\sm$, and
$P_a$ denotes the set of prevectors based at $a\in \tm$. We define a  \emph{prevector field}  on $\sm$ to be an \emph{internal} map $F:\sm \to \sm$ such that $(a,F(a)) \in P_a$ for every $a \in\tm$, that is, in coordinates 
 $ F(a) - a  \prec\DD $ for every $a \in \tm$.
If $F$ and $G$ are two prevector fields then we will say $F$ is equivalent to $G$ and write $F \equiv G$ if $F(a) \equiv G(a)$ for every $a\in\tm$ (recall Definition~\pr{dd3}).
A \emph{local} prevector field
 is an internal map $F:\s U \to \s V$ satisfying the above condition, where $U \su V \su M$ are open. When the distinction is needed, we will call a prevector field defined on all of $\sm$ a \emph{global} prevector field. 

The reason for allowing the values of a local  prevector field defined on $\s U$ to lie in a slightly larger range $\s V$ is in order to allow 
a prevector field $F$ to be restricted to a smaller domain which is not invariant under $F$. For example, if $M = \RR$ and one wants to restrict the prevector field $F$ given by $F(a)=a+\DD$, to the domain $\s (0,1)$, then one needs to allow a slightly larger range.
In the sequel we will usually not mention the larger range $V$ when describing a local prevector field, but it will always be tacitly assumed that we have such $V$ when needed. 
A second instance where it may be needed for the range to be slightly larger than the domain is the following natural setting for defining a
local  prevector field.

\begin{example}\label{e}
Let $p \in M$, $V$ a coordinate neighborhood of $p$ with image $V' \su \RR^n$, and $X$ a classical vector field on $V$, given in coordinates by $X' : V' \to \RR^n$. Then there is a neighborhood $U'$ of the image of $p$, with $U' \su V'$, such that we can define $F' : \s {U'} \to \s {V'}$ by 
$$F'(a)=a + \DD \cdot X'(a),$$ e.g. one can take  $U'$ such that $\overline{U'}$ is compact and $\overline{U'} \su V'$. 
For the corresponding $U \su V$ this induces a local prevector field $F: \s U \to \s V$ which realizes $X$ in
 $U$ in the sense of the following definition. (Recall Remark~\pr{67}.)
\end{example}

\begin{dfn}\label{realize}
A local prevector field $F$ on $\s U$ realizes the classical vector field $X$ on $U$ if
for every smooth $h:U \to \RR$  we have 
$Xh(a) = \vv{a \ F(a)}h$ for all $a\in U$. 
\end{dfn}

When realizing a vector field as in Example~\pr{e} it may indeed be necessary to restrict to a smaller neighborhood $U$, e.g. for $M=V=V'=(0,1)$ and $X' = 1$,  one needs to take $U=(0,r)$ for some
 $1 > r \in \RR$ in order for $F(a)=a+\DD$ to always lie in $\s V$. Note that Definition~\pr{realize} involves only  \emph{standard} points; see however Corollary~\pr{sr} for a discussion of this matter.

Different coordinates for the same neighborhood $U$ will induce equivalent realizations in $\s U$. More precisely, we show the following.

\begin{prop}\label{re}
For $U \su \RR^n$, let $X:U \to \RR^n$ be a classical vector field, $\f:U \to W \su \RR^n$ a change of coordinates, and $Y:W \to \RR^n$ the corresponding vector field, i.e. $Y(\f(a)) =  J_a X(a)$, where $J_a$ is the Jacobian matrix of $\f$ at $a$. Let $F,G$ be the prevector fields given by
$F(a) = a+\DD X(a)$, $G(a)=a+ \DD Y(a)$ as in Example~\pr{e}. Then $F \equiv \f^{-1} \circ G \circ \f$, or equivalently, 
$ \f \circ F(a) - G \circ\f(a)  \prec\prec \DD$ for all $a \in \tth U$.
\end{prop}

\begin{pf} 
Let $\f^i,Y^i,G^i$ be the $i$th component of $\f,Y,G$ respectively, and let $D^i_a$ be the differential of $\f^i$ at $a$, (so $D^i_a$ is the $i$th row of $J_a$).
Then we have
\begin{align*}
\f^i \circ F(a)  - G^i \circ\f(a)  &= \f^i(a+ \DD X(a)) -  \f^i(a) - \DD Y^i(\f(a))  \\  &= 
D^i_c\DD X(a) - \DD D^i_a X(a) = \DD( D^i_c - D^i_a) X(a) \prec\prec \DD,
\end{align*} 
where $D^i_c$ is the differential of $\f^i$ at some point $c$ in the interval between $a$ and $a+ \DD X(a)$. (Such $c$ exists by Remark~\pr{dh1}.)
Since this is true for each component $i$, we have $$ \f \circ F(a) - G \circ\f(a)  \prec\prec \DD.$$
\end{pf}

 \begin{dfn}\label{I} 
We define $I$ to be the \emph{identity} prevector field on $\sm$, (or on $\s U$ for any $U \su M$ or $U\su\rn$), i.e. $I(a)=a$ for all $a$.  
The prevector field $I$ corresponds to 
classical \emph{zero} vector field via the procedure of Example~\pr{e}.

\end{dfn}

\subsection{Regularity conditions}

If one wants to define various operations on prevector fields, such as their flow, or Lie bracket, then one must assume some regularity properties. 
Recall that a classical vector field $X:U \to \RR^n$ is called \emph{Lipschitz} if there is $K \in \RR$ such that $\| X(a)-X(b) \| \leq K\|a - b\|$ for $a,b \in U$.
For the local prevector field $F$ of Example~\pr{e}, where $F(a)-a=\DD X(a)$, this translates into $$\Bigl\| \Bigl( F(a)-a \Bigr) - \Bigl( F(b)-b \Bigr) \Bigr\| \leq K\DD\|a - b\|$$ for $a,b \in \s U$. 
This motivates the following definition.

\begin{dfn}\label{d1}
A prevector field $F$ on a smooth manifold $M$ is of class $D^1$ if  
whenever $a,b \in\tm$ and $(a,  b) \in P_a$ then in coordinates 
 $F(a)-a - F(b) + b \prec\DD \|a-b\| $. 
\end{dfn}

One can then also think of ``order $k$'' Lipschitz conditions on prevector fields, for the definition of which we will use the following Euler notation for finite differences. Given vector spaces $V,W$ (classical or nonstandard), and given $A \su V$, $a \in A$, and $v_1,\dots,v_k \in V$ such that $a+e_1v_1+\cdots+e_kv_k \in A$ for all
$(e_1,\dots,e_k)\in\{0,1\}^k$, and given a function $F:A \to W$, we define the $k$th difference 
$\xxk F(a)$ as follows:
$$\xxk F(a) 
= \sum_{(e_1,\dots,e_k)\in\{0,1\}^k}  (-1)^{\sum e_j}F(a+e_1v_1+\cdots+e_k v_k).$$
We note that in terms of this difference notation, the $D^1$ condition can be stated as follows:
$\x^1_{b-a}(F-I)(a)\prec\DD\|a-b\|$
 for any $a,b\in\tm$ with $a-b\prec\DD$, (recall that $I$ denotes the identity prevector field, i.e. $I(a)=a$ for all $a$).
Or, if we let $v=b-a$ then this can be written as $\x^1_v(F-I)(a)\prec\DD\|v\|$.

Generalizing to higher order differences, we define the $D^k$ regularity condition on a prevector field $F$ by the following condition, in coordinates in a neighborhood $U$:
For any $a \in \tu$ and any $v_1,\dots,v_k \in \s\RR^n$ with $v_i  \prec \DD$,
$$ \xxk (F-I)(a)    \prec\DD\|v_1\| \| v_2\|\cdots\|v_k\|.$$
We note that for $k \geq 2 $, $\xxk I(a)=0$, 
so for $k \geq 2$ the $D^k$ condition simplifies to:  $$\xxk F(a) \prec\DD\|v_1\| \| v_2\|\cdots\|v_k\|.$$
For $k=2$ this reads as follows.

\begin{dfn}\label{d2}
A prevector field $F$ on a smooth manifold $M$ is of class $D^2$   if for any $a \in \tm$, we have in coordinates that
for any $v,w \in \s\RR^n$ with $v,w  \prec \DD$,
$$\x^2_{v,w}F(a)=F(a) -F(a+v) - F(a+w) + F(a+v+w) \prec\DD\|v\|\|w\|.$$ 
\end{dfn}

We will show that the definitions of $D^1$ and $D^2$ prevector fields are independent of coordinates in Propositions 
\pr{5} and \pr{4} respectively.
We will  show in Proposition~\pr{3} that $D^2$ implies  $D^1$. 
In fact, the proof of the invariance of $D^2$ will use the fact  that $D^2$ implies $D^1$ in any given coordinates, which in turn relies on the technical Lemma~\pr{1}.
In Proposition~\pr{6} we will show that the prevector field of Example~\pr{e} induced  by a classical vector field of class $C^k$, is a $D^k$ prevector field. 
This is in fact the central motivation for our definition of $D^k$, but we note that our $D^k$ is in fact a weaker condition than $C^k$,
e.g. $F$ of Example~\pr{e} is $D^1$  if  $X$ is (order 1) Lipschitz, which is weaker than $C^1$.  
We remark that a definition of $D^0$ along the above lines would simply  amount to $F(a)-a \prec \DD$, i.e. $F$ being a prevector field. Note, however,  that in our definitions above, being a prevector field is part of the definition of $D^k$.

\begin{prop}\label{6}
For open $W \su \RR^n$, 
let  $X:W \to \RR^n$ be a classical $C^k$ vector field. 
Then for any $a \in \tth W$ and any $v_1,\dots,v_k \in \s\RR^n$ with $v_i  \prec\prec 1$, 
$$\xxk X(a)  \prec\|v_1\| \| v_2\|\cdots\|v_k\|.$$
It follows that if $F$ is the prevector field on $\s W$ of Example~\pr{e}, i.e. $F(a)-a=\DD X(a)$,
then $F$ is $D^k$. 
\end{prop}

\begin{pf}
Let $U \su W$ be a smaller neighborhood of $st(a)$ for which all $k$th partial derivatives of $X$ are bounded. 
Let $X^1,\dots,X^n$ be the components of $X$.
Given $p \in U$ and $v_1,\dots,v_k \in \RR^n$ such that $p+s_1v_1+\cdots+s_kv_k \in U$ 
for all $0 \leq s_1,\dots,s_k \leq 1$,   let $$\psi^i(s_1,\dots,s_k) = X^i(p+s_1v_1+\cdots+s_kv_k).$$ 
By iterating the mean value theorem $k$ times there is  $(t_1,\dots,t_k) \in [0,1]^k$ such that 
$$\sum_{(e_1,\dots,e_k)\in\{0,1\}^k}  (-1)^{\sum e_j}\psi^i(e_1,\dots,e_k)
={\pa^k \over \pa s_1 \pa s_2 \cdots \pa s_k}\psi^i(t_1,\dots,t_k)(-1)^k$$ 
(For the case $k=2$ see e.g. Rudin \ct[Theorem 9.40]{r}).
So 
\begin{align*}
|\xxk X^i(p)|&= 
\Bigl|\sum_{(e_1,\dots,e_k)\in\{0,1\}^k}  (-1)^{\sum e_j}X^i(p+e_1v_1+\cdots+e_kv_k) \Bigr| 
\\ &= \Bigl|\sum_{(e_1,\dots,e_k)\in\{0,1\}^k}  (-1)^{\sum e_j}\psi^i(e_1,\dots,e_k)\Bigr| 
\\&=\Bigl|{\pa^k \over \pa s_1\cdots \pa s_k}\psi^i(t_1,\dots,t_k)\Bigr| \leq C_i\|v_1\| \| v_2\|\cdots\|v_k\|
\end{align*}
where $C_i$ is determined by a bound for all $k$th partial derivatives of $X^i$ in $U$.

This is true for each $X^i$, $i=1,\dots,n$, and so there is a $K\in\RR$ such that 
$$\|\xxk X(p) \| \leq K  \|v_1\| \| v_2\|\cdots\|v_k\|$$ for every $p \in U$ 
and $v_1,\dots,v_k\in \RR^n$ such that $p+s_1v_1+\cdots+s_kv_k \in U$ 
for all $0 \leq s_1,\dots,s_k \leq 1$, $s_i \in \RR$.
By transfer the same is true, with the same $K$, 
for all $p \in \s U$ and $v_1,\dots,v_k \in\s\RR^n$ such that $a+s_1v_1+\cdots+s_kv_k \in \s U$
 for all  $0 \leq s_1,\dots,s_k \leq 1$, $s_i \in \s\RR$.
In particular this is true for our $a$  and all  $v_1,\dots,v_k \in \s\RR^n$ 
with $v_i  \prec\prec 1$.
\end{pf}

\begin{remark}\label{vfcc}
Proposition~\pr{6} was stated for a $C^k$ vector field $X:U \to \RR^n$, but for the first statement one can think of $X$ as any $C^k$ map, and 
indeed in the proof of  Proposition~\pr{4} below it will be used for $X=\f:U \to W$ a $C^k$ change of coordinates.
\end{remark}

We note that a $D^1$  prevector field $F$ satisfies the following.

\begin{prop}\label{0}
If $F$ is $D^1$ and $a,b \in\tm$ with $a-b \prec \DD$, then $$\|a-b\| \prec \|F(a)-F(b)\| \prec \|a-b\|.$$ 
More in detail, given $K\prec 1$ such that $ \|F(a)-F(b)-a+b \| \leq K\DD\|a-b\|$
we have $(1-K\DD)\|a-b\| \leq \|F(a)-F(b)\| \leq (1+K\DD)\|a-b\|$. 
\end{prop}

\begin{pf} We have $|\|F(a)-F(b)\|-\|a-b \| | \leq \|F(a)-F(b)-a+b \| \leq K\DD\|a-b\|$, so  
$(1-K\DD)\|a-b\| \leq \|F(a)-F(b)\| \leq (1+K\DD)\|a-b\|$. 
\end{pf}

\subsection{Invariance of regularity conditions}

The definitions of $D^1$ and $D^2$ above are in terms of coordinates. In the present section we prove that these definitions are in fact independent of coordinates, and that every $D^2$ prevector field is also $D^1$. This section is quite technical and can be skipped on first reading. Lemma \pr{1} that we prove and use in this section will be used again only in the proof of Lemma \pr{10}.

\begin{prop}\label{5}
The definition of $D^1$ is independent of coordinates.
\end{prop}

\begin{pf} 
Let $U,W \su \RR^n$ be two coordinate charts for a neighborhood of $st(a)$, and $\f:U \to W$ the change of coordinates map.  
Let $F$ be a $D^1$ prevector field in $U$ and $G$ the corresponding prevector field in $W$ i.e. 
$G \circ \f = \f \circ F$. For $b$ with $a-b \prec\DD$ we must show  $$G(\f(a)) - G(\f(b)) - \f(a) + \f(b) \prec \DD\|\f(a)-\f(b)\|.$$ 
Let $\phi = \f^i$ be the $i$th component of $\f$.
By Remark~\pr{dh1} 
there is a point $x$ on the interval between $F(a)$ and  $F(b)$ such that $\phi(F(a)) - \phi(F(b)) = D_x(F(a)-F(b))$, where $D$ is the differential of $\phi$.
There is a point $y$ on the interval between $a$ and $b$ such that $\phi(a) - \phi(b) = D_y(a-b)$.
So \begin{align*}\phi(F(a)) - \phi(F(b)) - \phi(a) + \phi(b) 
&= D_x(F(a)-F(b)) - D_y(a-b)
\\ &= D_x(F(a)-F(b)-a+b) - (D_y - D_x)(a-b) \\& \prec\DD\|a-b\| \prec \DD\|\f(a)-\f(b)\|,\end{align*} since 1) the entries of $D_x$ are 
finite, 2) $F(a)-F(b)-a+b \prec\DD\|a-b\|$ by assumption, 3)
$D_x-D_y \prec \|x-y\| \prec\DD$ (assuming the partial derivatives of $\varphi$ are Lipschitz, e.g. if $\varphi$ is $C^2$), and 4) the entries of the Jacobian of $\f^{-1}$ are finite, giving $\|a-b\| \prec \|\f(a)-\f(b)\|$. 

This is true for all components $\phi=\f^i$ of $\f$ and so it is true for $\f$, i.e.  
 $$\f(F(a)) - \f(F(b)) - \f(a) + \f(b)  \prec \DD\|\f(a)-\f(b)\|$$
which completes the proof since $\f \circ F =  G \circ \f $. 
\end{pf}

We would now like to show that for any given coordinates, $D^2$ implies $D^1$. 
We first prove the following technical lemma, which will also be used in the proof of Lemma~\pr{10}.
We demonstrate the content of this lemma with a simple example. Let $f,g:\s\RR \to\s\RR$ be $f(x)=cx$ and $g(x)=cd\sin {\pi  \over 2 d}x$, with $d\prec\prec 1$, then $f(0)=g(0)=0$. We now advance by steps of size $d$ and see how $f$ and $g$ develop. The increment of $f$ and $g$ after one step is the same, $f(d)=g(d)=cd$. But after $m$ steps with $m \geq {1 \over d}$, $f(md)\geq c$ whereas $g(md)\leq cd \prec\prec c$. The increments of $f$ properly accumulate along the $m$ steps to produce a large value for $f(md)$ in comparison to $f(d)$, due to the fact that the increment $f(x+d)-f(x)$ is constant. On the other hand, $g(md)$ remains small since the increments of $g$ are not sufficiently persistent, this being reflected in the fact that the difference $\Bigl(g(x+2d)-g(x+d)\Bigr) - \Bigl(g(x+d)-g(x)\Bigr)$ between successive increments
 is not sufficiently small compared to the first increment $g(d)-g(0)$.
This lemma will, in fact,  be used in the reverse direction, namely,  a bound on $f(md)-f(0)$ will be used 
in order to obtain a stronger bound on $f(d)-f(0)$.

\begin{lemma}\label{1}
Let $B \su \RR^n$ be an open ball around the origin $0$, 
let $a\in \ha(0)$, and let $0\neq v\in\s\RR^n$ with $v\prec\prec 1$. 
If $G:\s\! B \to \s \RR^n$ is an internal function satisfying
$$\x^2_{v,v}G(x) \prec \|v\| \|G(a) - G(a+v) \|$$ for all $x \in \tth \! B$, 
 then there is $m \in \s\N$ such that
$a + mv \in \tth\! B$ and  $$G(a) - G(a+v) \prec \|v \| \|G(a) - G(a+mv)\|.$$ 
\end{lemma}

\begin{pf}
Let $N=\lf r/\|v\| \rf$ where $0<r \in \RR$ is slightly smaller than the radius of $B$, and for $0 \leq x \in \s\RR$, $\lf x \rf \in\s\N$ is the integer part of $x$.
So any $m \leq N$ satisfies that $a+mv \in \tth \! B$. 
Let $A=G(a)-G(a+v)$. For $0 \leq j \leq N$ let $x_j=a+jv$, then by our assumption on $G$ we have
$G(x_j)-2G(x_{j+1})+G(x_{j+2})= \x^2_{v,v}G(x_j) =C_j \|v\| \|A \|$ with $C_j \prec 1$. Let $C$ be the maximum of $C_0,\cdots, C_N$, then $C\prec 1$ and $G(x_j)-2G(x_{j+1})+G(x_{j+2})  \leq C \|v\| \|A \|$ for all $0 \leq j \leq N$. Given  $k \leq N$ we have 
\begin{multline*} 
\Bigl\| A-\Bigl(G(x_k)  - G(x_{k+1})\Bigr) \Bigr\| = \Bigl\|\Bigl(G(x_0) - G(x_1)\Bigr)-\Bigl(G(x_k) - G(x_{k+1})\Bigr)\Bigr\| \\ \leq
\sum_{j=0}^{k-1}   \Bigl\|\Bigl(G(x_j) - G(x_{j+1})\Bigr) -  \Bigl(G(x_{j+1}) - G(x_{j+2})\Bigr)\Bigr\|    
\leq C k  \| v \| \|A\|. 
\end{multline*}
So for any $m \leq N$, 
$$\Bigl\|mA - \Bigl(G(x_0) - G(x_m)\Bigr) \Bigr\| = \Bigl\| \sum_{k=0}^{m-1} \Bigl( A-\Bigl(G(x_k) - G(x_{k+1})\Bigr)\Bigr)  \Bigr\|  \leq C m^2 \|v\|\|A\|,$$ so
 $\Bigl\|mA - \Bigl(G(x_0) - G(x_m)\Bigr) \Bigr\| = K m^2 \|v\|\|A\|$
with $K \prec 1$.
It follows that $$ m\|A\| - K m^2 \|v\|\|A\|   \leq  \|G(x_0) - G(x_m) \|,$$
and so, multiplying by $\|v\|$, we have
$ m\|v\| \|A\|(1 - K m \|v\|)  \leq \|v\| \|G(x_0) - G(x_m) \| $.

Now let $m = \min \{ \ N \ , \ \lf{1 \over  2K \|v\|}\rf  \ \}$, 
then 
$$ m\|v\|\|A\|/2     \leq \|v\| \|G(x_0) - G(x_m) \| .$$
By definition of $N$ and since $K\prec 1$ we 
have that $m\|v\|$ is appreciable, i.e. not infinitesimal, and so finally 
$A  \prec     \|v\| \|G(x_0) - G(x_m) \| $, that is, 
$G(a)-G(a+v)  \prec     \|v\| \|G(a) - G(a+mv) \|$.
\end{pf}

\begin{prop}\label{3}
If $F$ is  $D^2$ for some choice of coordinates in $W$, then $F$ is $D^1$.
\end{prop}

\begin{pf}
Given $a \in \tth W$, in the given coordinates take some ball $B$ around $st(a)$.
Define $G=F-I$, i.e. $G(x)=F(x)-x$, then we must show for any $v\prec\DD $ that $G(a)-G(a+v) \prec\DD \|v \|$ (here $v = b-a$ in Definition~\pr{d1}).  
If $G(a)-G(a+v)  \prec\prec \DD \|v \|$ then we are certainly done.
Otherwise $\DD v \prec\|G(a) - G(a+v)\|$ so $\DD \|v\|^2 \prec\|v\| \|G(a) - G(a+v)\|$,
and on the other hand $\x^2_{v,v}G(x)=\x^2_{v,v}F(x) \prec\DD\|v\|^2$ for all $x$, since $\x^2_{v,v}I(x)=0$, and by taking $v=w$ in Definition~\pr{d2}.
Together we have 
$\x^2_{v,v}G(x) \prec \|v \|\|G(a) - G(a+v)\|$,
so by Lemma~\pr{1} there is $m\in\s\N$ such that $a+mv \in \tth\! B$ and 
$$G(a) - G(a+v) \prec\|v\| \|G(a) - G(a+mv) \|
\leq \|v\| \Bigl(\|G(a)\| + \|G(a+mv)\|\Bigr) \prec \|v\|\DD$$ since $F$ is a prevector field 
and so $G(x)=F(x)-x\prec\DD$ for all $x$.
\end{pf}

We need  one more lemma before proving that $D^2$ is independent of coordinates.

\begin{lemma}\label{2}
In given coordinates, $F$ is  $D^2$     if and only if  it satisfies 
$$\xt F(a)  \prec\DD \bigl(\max \{\|v\|,\|w\|\}\bigr)^2$$
for every $a$, and every $v,w \prec\DD$.
\end{lemma}

\begin{pf}
Clearly $D^2$ implies the above condition. For the converse, 
say $\|v \| \leq \|w \|$. If $w \prec\|v\|$ then the two conditions are clearly equivalent. 
Otherwise let $n = \lf\|w\| / \| v \|\rf$, and $w'=w/n$. Then $\| v \| \leq \|w' \| \leq {n+1 \over n} \|v \|$, and so by the preceding remark the two conditions are equivalent for $v,w'$, and so for every $1 \leq k \leq n$,
$$
C_k = {\x^2_{v,w'}F(a+(k-1) w') \over \DD\|v\|\|w'\|}
$$
is finite. Let $C$ be the maximum of $C_1,\dots, C_n$ then $C$ is finite and 
$$ \x^2_{v,w'}F(a+(k-1) w')  \leq C\DD\|v\|\|w'\|$$
for every $1 \leq k \leq n$. And so 
\begin{align*}
\| \xt & F(a)\|= 
\|F(a) -F(a+v) - F(a+w) + F(a+v+w) \| \\ \leq &\sum_{k=1}^n
\| F(a+(k-1) w') -F(a+(k-1) w'+v) - F(a+k w') + F(a+k w'+v) \|
\\ =& \sum_{k=1}^n  \| \x^2_{v,w'}F(a+(k-1) w')   \| \leq nC\DD\|v\|\|w'\|  = C\DD\|v\|\|w\|.
\end{align*}
\end{pf}

\begin{prop}\label{4}
The definition of $D^2$ is independent of coordinates.
\end{prop}

\begin{pf}
Let $\f:U \to W$ be a change of coordinates.  Let $F$ be a $D^2$ prevector field in $U$ and $G$ the corresponding prevector field in $W$ i.e. 
$G \circ \f = \f \circ F$.
Given $p \in {}^\ha W$, and $x,y \in\s\RR^n$ with $x,y\prec\DD$, and say $ \|x \|\leq \|y\|$, then by Lemma~\pr{2}    
it is enough to show $\x^2_{x,y} G(p)=G(p)-G(p+x)-G(p+y)+G(p+x+y) \prec\DD \|y\|^2$.
Let $a,v,w$ be such that $\f(a)=p$, $\f(a+v)=p+x$, $\f(a+w)=p+y$.
Then  $\|w\|\prec \|y\|$ and so it is enough to show 
\begin{equation*}\label{x1}\tag{1}
\f( F(a))-\f(F(a+v))-\f(F(a+w))+G(p+x+y) \prec\DD \|w\|^2.  
\end{equation*}  
Since $F$ is $D^2$, by Proposition~\pr{3} it is also $D^1$, and so by Proposition~\pr{5} $G$ is $D^1$. 
So \begin{align*}\label{x2}\tag{2}
G(p+x&+y)   -  G( \f(a+v+w))-(p+x+y) + \f(a+v+w)  \\ & \prec\DD\| p+x+y-\f(a+v+w)\| \\ & = \DD\| -\f(a)+\f(a+v)+\f(a+w)-\f(a+v+w)\|  \prec\DD\|w\|^2,   \end{align*} 
by Remark~\pr{vfcc} (assuming $\f$ is $C^2$) and  since
$p+x+y=-\f(a)+\f(a+v)+\f(a+w)$ and  $\|v\|  \prec \| w\|$.
In view of (\pr{x2}) we see that (\pr{x1}) holds     
if and only if  
\begin{multline*}
\f( F(a))-\f( F(a+v))-\f( F(a+w)) 
+\f (  F (a+v+w))   \\  -\f(a)+\f(a+v)+\f(a+w) - \f(a+v+w)  
\\ =  \xt(\f\circ F -\f)(a)    \prec\DD \|w\|^2, 
\end{multline*}
so we proceed to prove this last inequality. 
Let $\phi = \f^i$ be the $i$th component of $\f$. We have
\begin{multline*}\label{x3}\tag{3}
\phi( F(a+v)) - \phi( F(a)) = \\ D(F(a+v)-F(a))+    {1\over 2}(F(a+v)-F(a))^tH_1(F(a+v)-F(a))
\end{multline*}
where $D=D_{F(a)}$ is the differential of $\phi$ at $F(a)$ and $H_1$ is the Hessian matrix of $\phi$ at some point on the interval between $F(a)$ and $F(a+v)$, (recall Remark~\pr{dh2}). Similarly 
\begin{multline*}\label{x4}\tag{4}
\phi( F(a+w)) - \phi( F(a)) =\\ D(F(a+w)-F(a))+    {1\over 2}(F(a+w)-F(a))^tH_2(F(a+w)-F(a))
\end{multline*}
and 
\begin{multline*}\label{x5}\tag{5}
\phi( F(a+v+w)) - \phi( F(a)) =\\ D(F(a+v+w)-F(a))+    {1\over 2}(F(a+v+w)-F(a))^tH_3(F(a+v+w)-F(a))\end{multline*}
with $H_2,H_3$ similarly defined. We have 
$$\phi(a+v) - \phi(a) = D_a v+    {1\over 2}v^tH'_1v$$ 
where $D_a$ is the differential of $\phi$ at $a$ and 
$H'_1$ is  the Hessian matrix of $\phi$ at some point on the interval between $a$ and $a+v$. Now let $Y_1 = H'_1 -H_1$ then $Y_1 \prec \DD$ (assuming the second partial derivatives of $\f$ are Lipschitz, e.g. if $\f$ is $C^3$), and we have 
\begin{equation*}\label{x6}\tag{6}
\phi(a+v) - \phi(a) = D_av+   {1\over 2} v^t(H_1+Y_1)v.
\end{equation*}
Similarly there are $Y_2,Y_3 \prec \DD$
such that 
\begin{equation*}\label{x7}\tag{7}
\phi(a+w) - \phi(a) = D_a w+   {1\over 2} w^t(H_2+Y_2)w
\end{equation*}
and
\begin{equation*}\label{x8}\tag{8}
\phi(a+v+w) - \phi(a) = D_a (v+w)+    {1\over 2}(v+w)^t(H_3+Y_3)(v+w)
\end{equation*}
Furthermore, since by Proposition~\pr{3} the prevector field $F$ is $D^1$,  we have $$F(a)-F(a+v)-a+(a+v) \prec\DD\| v\|$$ i.e. 
$F(a+v)-F(a)=v+\delta_1$ with $\delta_1 \prec \DD \|v\| \prec \DD \|w\|$. 
Similarly there are $\delta_2,\delta_3 \prec \DD \| w \|$ such that 
$F(a+w)-F(a)=w+\delta_2$, $F(a+v+w)-F(a)=v+w+\delta_3$. 
Substituting this  into the quadratic terms of  (\pr{x3}),(\pr{x4}),(\pr{x5}) we get

\begin{equation*}\label{x9}\tag{9}
\phi( F(a+v)) - \phi( F(a)) =  D(F(a+v)-F(a))+    {1\over 2}(v+\delta_1)^tH_1(v+\delta_1)
\end{equation*}
\begin{equation*}\label{x10}\tag{10}
\phi( F(a+w)) - \phi( F(a)) = D(F(a+w)-F(a))+    {1\over 2}(w+\delta_2)^tH_2(w+\delta_2)
\end{equation*}
\begin{equation*}\label{x11}\tag{11}
\phi( F(a+v+w)) - \phi( F(a)) = D(F(a+v+w)-F(a))+{1\over 2}(v+w+\delta_3)^tH_3(v+w+\delta_3)
\end{equation*}

Now
\begin{align*} 
 &\xt(\phi\circ F -\phi)(a)  = \\ & 
\phi( F(a))-\phi(F(a+v))-\phi(F(a+w))   +\phi (F (a+v+w)) \\&
  -\phi(a)+\phi(a+v)+\phi(a+w) - \phi(a+v+w) \\
=&  -  \Bigl(\phi( F(a+v)) - \phi( F(a))\Bigr)   
-  \Bigl(\phi( F(a+w)) - \phi(F(a))\Bigr)  
+  \Bigl(\phi( F(a+v+w)) - \phi( F(a))\Bigr)  \\&
+\Bigl(\phi(a+v) - \phi(a)\Bigr) 
+\Bigl(\phi(a+w) - \phi(a)\Bigr) 
-\Bigl(\phi(a+v+w) - \phi(a)\Bigr)
\end{align*}

Substituting (\pr{x6}),(\pr{x7}),(\pr{x8}),(\pr{x9}),(\pr{x10}),(\pr{x11}) for these six parenthesized summands, after all cancellations we remain with
\begin{multline*}
D\Bigl(F(a)-F(a+v)-F(a+w)+F(a+v+w)\Bigr) \\ 
- v^tH_1\delta_1  - w^tH_2\delta_2  + (v+w)^tH_3\delta_3  - {1\over 2}\delta_1^tH_1\delta_1 - {1\over 2}\delta_2^tH_2\delta_2  + {1\over 2}\delta_3^tH_3\delta_3 \\
+{1\over 2}v^tY_1v + {1\over 2}w^tY_2w - {1\over 2}(v+w)^tY_3(v+w) \prec \DD\|w\|^2.
\end{multline*}

This is true for all components $\phi=\f^i$ of $\f$ and so it is true for $\f$.
\end{pf}

\subsection{Operations on prevector fields}

We now show that addition of the vectors corresponding to $D^1$  prevector fields $F,G$ is realized by their composition $F \circ G$. More precisely, we show the following.

\begin{prop}
Let $F$ be a $D^1$ prevector field and $G$ any prevector field. Then for every $a\in\tm$, 
$$\vv{a \ F(G(a))} = \vv{a  \ F(a)}+\vv{a  \ G(a)}.$$
In particular if both $F,G$ are  $D^1$ then $F\circ G \equiv G \circ F$.
\end{prop}

\begin{pf}
In coordinates $\vv{a \ F(a)}+\vv{a \ G(a)} = \vv{a \ x}$ where $x=a+ (F(a)-a) + (G(a)-a) = F(a)+G(a)-a$.
So $F(G(a)) - x = F(G(a))-F(a) - G(a)+a \prec \DD\|G(a)-a\| \prec\prec \DD$. 
\end{pf}

We next show that  the composition of  $D^1$ (resp. $D^2$) prevector fields is  $D^1$ (resp. $D^2$).

\begin{prop}\label{40}
If $F,G$ are  $D^1$ then $F\circ G$ is $D^1$.
\end{prop}

\begin{pf} We have 
\begin{align*}\|F\circ G(a) &-F\circ G(b)-a+b\| \\&\leq \|F\circ G(a)-F\circ G(b) -G(a) + G(b)\| +\|G(a) -G(b) -a+b\| \\&
\prec \DD\|G(a)-G(b)\| + \DD\|a-b\| \prec \DD\|a-b\| \end{align*}  by Proposition~\pr{0}. 
\end{pf}

\begin{prop}\label{q}
If $F,G$ are  $D^2$ then $F\circ G$ is $D^2$.
\end{prop}

\begin{pf} In some coordinates 
let $p=G(a)$, $x=G(a+v)-G(a)$ and $y=G(a+w)-G(a)$, and so by Propositions \pr{3} and \pr{0} we have
$x \prec \|v\|$ and $y \prec \|w\|$. Also by Propositions \pr{3} and \pr{0} we have   
\begin{multline*} \|F(p+x+y)-F(G(a+v+w))\| \prec \| p+x+y-G(a+v+w) \| \\ = \|-G(a)+G(a+v)+G(a+w)-G(a+v+w)\| 
\prec \DD\|v\|\|w\|.\end{multline*} 
Now \begin{align*}
\|&\xt (F \circ G)(a)\| \\&=
\|F \circ G (a) -   F \circ G (a+v)  -         F \circ G ( a+w)   +              F \circ G (a+v+w) \|  \\ &= 
\| F (p) -   F (p+x)  -         F  ( p+y)   +              F \circ G (a+v+w) \|  \\
&\leq \| F (p) -   F (p+x)  -         F  ( p+y)   +  F(p+x+y)   \|  +  \|F(p+x+y)-       F (G (a+v+w))  \|  \\
&\prec \DD\|x\| \|y\| + \DD\|v\|\|w\| \prec \DD\|v\|\|w\|.\end{align*}
\end{pf}

\section{Global properties of prevector fields}\label{global}

The definition of $D^1$ and $D^2$ prevector fields relates to points in $\tm$ which are infinitely close to each other, in fact of distance $\prec\DD$. In this section we establish properties of $D^1$ and $D^2$ prevector fields valid on  appreciable neighborhoods, or on the whole of $\tm$.

\begin{prop}\label{18}
Let $F$ be a prevector field. 
\begin{enumerate}
\item If $W$ is a coordinate neighborhood with image $U \su \RR^n$  and $B \su U$ is a closed ball, 
then there is a finite  $C$ such that $\|F(a)-a\| \leq C\DD$ for all $a\in \s\! B$. (We use $F$ to denote both the prevector field itself, and its action in coordinates.)
\item If $G$ is another prevector field, then there is a finite $\beta$ such that   $\|F(a)-G(a)\| \leq \beta\DD$ for all $a\in \s\! B$.
\item If furthermore $F\equiv G$ then an infinitesimal such $\beta$ exists. 
\end{enumerate}
\end{prop}

\begin{pf}
The first statement is a special case of the second, by taking $G(a)=a$ for all $a$. So we prove the second statement.
Let $$A= \{ n \in \s\N :  \| F(a)-G(a)  \| \leq n\DD   \hbox{\ \ for every\ \ } a \in \s\! B \}.$$ Every infinite $n\in\s\N$ is in $A$ and so by underspill 
\footnote{Recall that for $\s\N$, ``underspill'' is 
the fact that 
if $A \su \s\N$ is an internal set, and $A$ contains all infinite $n$ then it must also contain a finite $n$.}
there is a finite $C$ in $A$.  

For $F \equiv G$, let 
$$A= \{ n \in \s\N :  \| F(a)-G(a)  \| \leq {\DD\over n}   \hbox{\ \ for every\ \ } a \in \s\! B \}.$$ Every
finite $n \in\s\N$ is in $A$ and so
by overspill \footnote{If $B \su \s\N$ is internal, and $B \supseteq \N$, then there must also be an infinite $n \in B$.} there is an infinite $n\in\s\N$ in $A$, and take $\beta = {1 \over n} \prec\prec 1$.
\end{pf}

\begin{prop}\label{20}
Let $F$ be a  $D^1$ prevector field. If $W$ is a coordinate neighborhood with image $U \su \RR^n$  and $B \su U$ is a closed ball, 
then there is $K \in \RR$ such that
\begin{equation*}
\| F(a)-a-F(b)+b\| \leq K\DD\|a-b\|  \tag{1}
\end{equation*} 
for all $a,b \in \s \! B$. 

It follows that
\begin{equation*}
(1-K\DD)\|a-b\| \leq \|F(a)-F(b)\| \leq (1+K\DD)\|a-b\|.   \tag{2}
\end{equation*}
\end{prop}

\begin{pf}

Let $N = \lf 1 / \DD\rf$. 
Given $a,b \in \s \! B$, for $k=0,\dots,N$ let
$a_k=a+{k \over N}(b-a)$, then $a_k - a_{k+1} \prec\DD$. Let $C_{ab}$ be the maximum of 
$${\| F(a_k)-a_k - F(a_{k+1}) + a_{k+1} \| \over \DD\|a_k-a_{k+1}\|}$$ for $k=0,\dots,N-1$,
then $C_{ab}$ is finite. 
For every $0\leq k \leq N-1$ 
we have  
 $$\| F(a_k)-a_k - F(a_{k+1}) + a_{k+1} \| \leq C_{ab}\DD\|a_k-a_{k+1}\|=C_{ab}\DD{\|a-b\|\over N}.$$
So $$\| F(a)-a - F(b) + b \| \leq \sum_{k=0}^{N-1} \| F(a_k)-a_k - F(a_{k+1}) + a_{k+1} \|  \leq C_{ab} \DD\|a-b\|.$$
Now let $$A=\{ n \in \s\N : \| F(a)-a - F(b) + b \| \leq n \DD\|a-b\| \hbox{\ \ for every\ \ } a,b \in \s\! B \}.$$ 
Since each $C_{ab}$ is finite, every infinite $ n \in \s\N$ is in $A$, and so by underspill, there is a finite $K$ in $A$,
and the first statement follows. 
The second statement follows from the first as in the proof of Proposition~\pr{0}.
\end{pf}

\begin{prop}\label{202}
Let $F$ be a  $D^2$ prevector field. If $W$ is a coordinate neighborhood with image $U \su \RR^n$  and $B \su U$ is a closed ball, 
then there is $K \in \RR$ such that
$$\| \xt F(a) \| \leq K\DD \|v\| \|w\|$$
 for all $a \in \s\! B$ and $v,w \in\s \RR^n$ such that $a+v,a+w,a+v+w\in\s\! B$. 
\end{prop}

\begin{pf}
The proof is similar to that of Proposition~\pr{20}.
Let $N = \lf 1 / \DD\rf$. 
Given $a \in\s\! B$ and $v,w\in\s\RR^n$ such that  $a+v,a+w,a+v+w\in\s\! B$, let 
$a_{k,l}=a+{k \over N}v + {l \over N}w$, $0 \leq k,l \leq N$.
Let $C_{avw}$ be the maximum of 
$${\|  \x^2_{v/N,w/N}F(a_{k,l})   \| \over \DD\|v/N\|\|w/N\|}$$ for $0\leq k,l \leq N-1$,
then $C_{avw}$ is finite. 
For every $0 \leq k,l \leq N-1$ we have   
\begin{multline*}   \| F(a_{k,l})-F(a_{k+1,l}) - F(a_{k,l+1}) + F(a_{k+1,l+1} \| = \\ \|  \x^2_{v/N,w/N}F(a_{k,l})   \|\leq C_{avw}\DD\|v/N\|\|w/N\|.\end{multline*}
Summing over  $0 \leq k,l \leq N-1$ we get 
$$\| F(a)-F(a+v)-F(a+w)+F(a+v+w)\| \leq C_{avw}\DD \|v\| \|w\|.$$
By underspill as in the proof of Proposition~\pr{20},  there is a single finite  $K$ which works for all $a,v,w$.
\end{pf}

When speaking about local $D^1$ or $D^2$ prevector fields, whenever needed we will assume, perhaps by passing to a smaller domain, that a constant $K$ as in Propositions \pr{20}, \pr{202} exists.

\begin{cor}\label{inj}
If $F$ is a  $D^1$ prevector field then $F$ is injective on $\tm$.
\end{cor}
 
\begin{pf}
Let $a \neq b \in \tm$. If $st(a)\neq st(b)$ then clearly $F(a) \neq F(b)$. Otherwise there exists a $B$ containing $a,b$ as in Proposition~\pr{20}, 
and (2) of that proposition implies $F(a) \neq F(b)$.
\end{pf}

We will now show that a $D^1$ prevector field is in fact bijective on $\tm$. We first prove local surjectivity, as follows.

\begin{prop}\label{surj}
Let $B_1 \su B_2 \su B_3 \su \RR^n$ be closed balls centered at the origin, of radii $r_1<r_2<r_3$. If $F:\s\! B_2 \to \s\! B_3$ is a local $D^1$ prevector field then $F(\s\! B_2) \supseteq \s\! B_1$.  
\end{prop}

\begin{pf}
Fix $0<s\in\RR$ smaller than $r_2-r_1$ and $r_3-r_2$. We will apply transfer to the following fact: For every function $f:B_2 \to B_3$, if 
$\|f(x)-f(y)\| \leq 2\|x-y\|$ for all $x,y \in B_2$ and $\|f(x)-x\| < s$ for all $x \in B_2$ then $f(B_2) \supseteq B_1$. This fact is indeed true since our assumptions on $f$ imply that it is continuous, and that for every $x \in \pa B_2$ the straight interval between $x$ and $f(x)$ is included in $B_3 -B_1$,
so $f|_{\pa B_2}$ is homotopic  in $B_3-B_1$ to the inclusion of $\pa B_2$. Now if some $p \in B_1$ is not in $f(B_2)$ 
then $f|_{\pa B_2}$ is null-homotopic in $B_3-\{ p\}$, 
and so the same is true for the inclusion of $\pa B_2$, a contradiction.
Applying transfer we get that for every internal function $f:\s\! B_2 \to \s\! B_3$,  if 
$\|f(x)-f(y)\| \leq 2\|x-y\|$  for all $x,y \in \s\! B_2$ and $\|f(x)-x\| < s$ for all $x \in \s\! B_2$ 
then $f(\s\! B_2) \supseteq \s\! B_1$.
In particular this is true for a $D^1$ prevector field $F:\s\! B_2 \to \s\! B_3$, by Proposition~\pr{20}(2).
\end{pf}

The following is immediate from Corollary~\pr{inj} and Propostion~\pr{surj}.

\begin{thm}\label{bij}
If $F: \sm \to \sm$ is a $D^1$ prevector field then $F|_{\tm} :\tm \to \tm$ is bijective.
\end{thm}

\begin{remark}\label{finv1}
On all of $\sm$, a $D^1$ prevector field may be noninjective and nonsurjective, e.g. take $M=(0,1)$ and  $F:\sm \to \sm$ given by $F(x)=\DD$ for $x \leq \DD$ and $F(x)=x$ otherwise. (Recall that the definition of $D^1$ prevector field imposes no restrictions at points of $\sm - \tm$).
\end{remark}

\begin{remark}\label{finv2}
 For $D^1$ prevector field $F$, the map $F|_\tm :\tm \to \tm$ and its  inverse  $(F|_{\tm})^{-1}$ are not internal if $M$ is noncompact, since their domain is not internal. 
On the other hand, for any $A \su M$, $F|_{\s\! A}$ is internal. Furthermore, on $\s\! B_1$ of Proposition~\pr{surj}, $F$ has an inverse $F^{-1}:\s\! B_1 \to \s\! B_2$ in the sense that  $F\circ F^{-1}(a) = a$ for all $a\in \s\! B_1$, and $F^{-1}$ is internal. So, for a local $D^1$ prevector field $F:\s U \to \s V$ we may always assume (perhaps for slightly smaller domain) that $F^{-1} : \s U \to \s V$ also exists, in the above sense. As mentioned, we will usually not mention the range $\s V$ but rather speak of a local prevector field on $\s U$.
\end{remark}

\begin{prop}\label{30}
If $F$ is  $D^1$ then $F^{-1}$ is $D^1$. ($F^{-1}$ exists by Remark~\pr{finv2}.) 

More in detail, if for  $x=F^{-1}(a), y=F^{-1}(b)$ there is given $K \prec 1$ such that   
$\|F(x)-F(y)-x+y\| \leq K\DD\|x-y\|$, then $\|F^{-1}(a)-F^{-1}(b)-a+b\| \leq K'\DD\|a-b\|$, with $K'$ only slightly larger, namely
$K' = K/(1-K\DD)$.
\end{prop}

\begin{pf}
Let $x=F^{-1}(a)$,  $y=F^{-1}(b)$, then 
\begin{multline*}\|F^{-1}(a)-F^{-1}(b)-a+b\| =\\ \|x-y-F(x)+F(y)\| \leq K \DD\|x-y\| \leq K' \DD\|F(x)-F(y)\|=K'\DD\|a-b\|\end{multline*}
by Lemma~\pr{0}. 
\end{pf}

We conclude this section with the following observations.

\begin{lemma}\label{103}
Let $F,G$ be $D^1$ prevector fields. If   $F(a) \equiv G(a)$ for all \emph{standard} $a$, 
then $F(b)\equiv G(b)$ for all nearstandard $b$, i.e. $F\equiv G$.
\end{lemma}

\begin{pf}
Given $b$, 
let $K$ be as in Proposition~\pr{20} for both $F$ and $G$, in a ball around $a=st(b)$. Then 
\begin{align*}\|F(b)-G(b)\| &= \| F(b)-F(a)-b+a+F(a) - G(a) +G(a) -G(b) -a +b \|  \\&
\leq \|F(b)-F(a)-b+a\| + \|F(a)-G(a)\| + \|G(a)-G(b)-a+b\|  \\ & \leq K\DD\|a-b\| +\|F(a)-G(a)\| + K\DD\|a-b\|  \prec\prec\DD. \end{align*}
\end{pf}

Recall that Definition~\pr{realize}, which defines when a prevector field $F$ realizes a classical vector field $X$,  involves only  \emph{standard} points.
It follows from Lemma~\pr{103} that if $F$ is $D^1$ then this determines $F$ up to equivalence.
Namely, we have the following.

\begin{cor}\label{sr}
Let $U \su\rn$ be open, and $X:U\to\rn$ a classical vector field. 
If $F,G$ are two $D^1$ prevector fields that realize $X$ then $F\equiv G$. 
In particular, if $X$ is Lipschitz and $G$ is a $D^1$ prevector field that realizes $X$, then
$F\equiv G$, where $F$ 
is the prevector field obtained from 
$X$ as in Example~\pr{e}.
\end{cor}

\section{The flow of a prevector field}\label{f}

In this section we define and study the flow of global and local prevector fields. 
In our definition of a prevector field as a map from $\sm$ to itself, we wish to view $F$ as its
own flow at time $\DD$. The flow for later time $t$ should thus be defined by iterating $F$ the appropriate number of times. (Thus the classical notion of a vector field being the infinitesimal generator of its flow receives literal meaning in our setting.)  

Thus, for a global prevector field $F: \sm \to \sm$ and for $0 \leq t \in\s\RR$, let $n = n(t) = \lf t/\DD \rf$ and define the flow $F_t$ of $F$ at time 
$t$ to be $F_t(a)=F^n(a)$, where $F^n$ is given by the map 
$\s{\mbox{Map}} (M) \times \s \N \to \s{\mbox{Map}}(M)$  
which is the extension of the map $\mbox{Map}(M) \times \N \to \mbox{Map}(M)$ 
taking $(f,n)$ to  $f^n=f \circ f\circ \cdots \circ f$. 
The flow of a local prevector field is similarly defined, only a bit of care is needed regarding its domain. 
So, for local prevector field $F:\s U \to \s V$,  extend $F$ to $F':\s V \to \s V$ by defining $F'(a)=a$ for all $a \in \s V - \s U$, 
and let 
$Y_n = \{ a \in \s U \ : \  (F')^n(a) \in \s U  \}$. We set the domain of $F_t$ to be $Y_{n(t)}$, 
where it is defined by $F_t(a) = F'_t(a)$. 

We would also like to consider $F_t$ for $t \leq 0$. For global prevector field $F$ which is bijective on $\sm$, or for local prevector field which is bijective in the sense of Remark~\pr{finv2}, in particular a $D^1$ local prevector field, we define $F_t$ for $t \leq 0$
to be $(F^{-1})_{-t}$. 

Note that for any \emph{global} prevector field $F$,  $F_t$ is defined for all $t \geq 0$, unlike the situation for the classical flow of a classical vector field. Similarly $F_t$ is defined for all $t \leq 0$ if $F$ is bijective.

Directly from the definition of a flow, we may immediately notice the following.

\begin{prop}\label{invfl}
A prevector field $F$ is invariant under its own flow $F_t$, where the action of a map $h$ on a prevector $(a,x)$ is given, as in Definition~\pr{dd5}, by $(h(a),h(x))$.
\end{prop}

\begin{pf}
Let $n = \lf t/\DD \rf$ then $F_t((a,F(a))) = F^n((a,F(a)))=(F^n(a),F^{n+1}(a)) = (b,F(b))$ where $b=F^n(a)=F_t(a)$. 
\end{pf}

 \subsection{Dependence on initial condition and on  prevector field}\label{sf1}

We now establish bounds on the distance in coordinates between two flows $F_t(a),F_t(b)$ of a given $D^1$ prevector field $F$, and between the flows $F_t(a),G_t(a)$ of two different prevector fields. These bounds  can of course be combined into a bound on the distance between $F_t(a)$ and $G_t(b)$.

\begin{thm}\label{flow1}
Let $F$ be a local $D^1$ prevector field on $\s U$.
Given $p \in U$ and a coordinate neighborhood of $p$ with image $W \su \RR^n$, 
let $B' \su B \su W$ be closed balls of radii $r/2,r$ around the image of $p$.  Suppose 
$\|F(a)-a-F(b)+b\| \leq K\DD\|a-b\|$ for all $a,b \in \s\! B$, with $K$ a  finite constant (such finite $K$ exists by Proposition~\pr{20}),
then there is
  $0 < T \in \RR$ 
such that $F_t(a) \in \s\! B$ for all $a \in \s\! B'$ and $-T \leq t \leq T$. Furthermore, 
for all $a,b \in \s\! B'$ and $0 \leq t \leq T$:
$\|F_t(a) - F_t(b) \| \leq e^{Kt}\|a-b\|$.

If we take a slightly larger constant $K'=K/(1-K\DD)^2$, then
for all $a,b \in \s\! B'$ and $-T \leq t \leq T$:
$$e^{-K'|t|}\|a-b\|\leq \|F_t(a) - F_t(b) \| \leq e^{K'|t|}\|a-b\|.$$
\end{thm}

\begin{pf} 
We first prove the statement for $t\geq 0$. 
Let $C$ be as in Proposition~\pr{18}(1).
Take  $T = {r \over 2C}$, then $0 < T \in \RR$, and we have 
by internal induction \footnote{If $A$ is an internal subset of $\s\N$ that contains $1$
and is closed under the successor function $n \mapsto n + 1$, then $A = \s\N$.
So, one can prove by induction in $\s\N$, as long as all objects under discussion are internal. This is called \emph{internal induction}.}
for $a \in \s\! B'$, $0 \leq t \leq T$, and $n=\lf t / \DD \rf$, that $\|F^n(a) - a\| \leq \sum_{m=1}^n \|F^m(a) - F^{m-1}(a)\| \leq nC\DD \leq {r \over 2}$, and so $F^n(a) \in \s\! B$. 
By Proposition~\pr{20}(2) we have 
$(1-K\DD)\|a-b\| \leq \|F(a)-F(b)\| \leq (1+K\DD)\|a-b\|$ for all $a,b \in \s\! B$, and so by internal induction 
$$(1-K\DD)^n\|a-b\|\leq \|F^n(a) - F^n(b) \| \leq (1+K\DD)^n\|a-b\|.$$
For $t\geq 0$ we have $(1+K\DD)^n \leq e^{Kt}$ since $1+K\DD \leq e^{K\DD}$, 
and letting $h=\frac{1}{1-K\DD}$ we have $e^{-hKt}\leq (1-K\DD)^n$
since $e^{-hK\DD} \leq 1-hK\DD + \frac{h^2 K^2\DD^2}{2}
 \leq 1-hK\DD +hK^2\DD^2=1-K\DD$. For $t\leq 0$ we are considering $F^{-1}$. The same constant $C$ can be used, and by Proposition~\pr{30} $K$ should be replace by $hK$, and so $hK$ is replaced by $K'=h^2 K$.
\end{pf}

\begin{thm}\label{flow2}
Let $F$  be a $D^1$ local prevector field on $\s U$ and let $G$ be any local prevector field on $\s U$. 
Given a coordinate  neighborhood included in $U$ with image $W \su\RR^n$, 
let $A'\su A \su \s W$ be internal sets. 
Suppose
\begin{enumerate}
\item
$\|F(a)-G(a)\| \leq \beta\DD$ for all $a\in  A$, with some constant  $\beta$. (If $F\equiv G$ then 
an infinitesimal such $\beta$ exists by Proposition~\pr{18}(3)),
\item  
$\| F(a)-a-F(b)+b\| \leq K\DD\|a-b\|$
for all $a,b \in A$, with $K$ a finite constant. 
(Such finite $K$ exists by Proposition~\pr{20}), 
\item $0<T\in\RR$ is such that $F_t(a)$ and $G_t(a)$ are in $ A$ 
for all $a \in  A'$ and $0 \leq t \leq T$.
\end{enumerate}
Then for all
$a\in  A'$ and $0 \leq t \leq T$,
$$\|F_t(a) - G_t(a) \| \leq {\beta\over K} (e^{Kt}-1)\leq \beta t e^{Kt}.$$

If $G^{-1}$ exists, e.g. if $G$ is also $D^1$, and if 
$F_t(a)$ and $G_t(a)$ are in $A$ 
for all $a \in  A'$ and $-T \leq t \leq T$,
then
$\|F_t(a) - G_t(a) \| \leq {\beta\over K} (e^{K|t|}-1) \leq \beta |t| e^{K|t|}$
for all $-T \leq t \leq T$. 
\end{thm}

\begin{pf} 
Again it is enough to prove the statement  for positive $t$. We prove by internal induction that
$$ \|F^n(a) - G^n(a) \| \leq {\beta\over K} \Bigl((1+K\DD)^n - 1\Bigr) $$
which implies the statement.
By  Proposition~\pr{20}(2) we have 
\begin{align*}  \|F^{n+1}(a)-G^{n+1}(a)   \| &\leq   \|   F(F^n(a))-F(G^n(a)) \|  + \|   F(G^n(a))-G(G^n(a)) \|  \\
& \leq (1+K\DD) \|F^n(a)-G^n(a)\| + \beta\DD \end{align*}
from which the induction step from $n$ to $n+1$ follows.
\end{pf}

\begin{cor}\label{flow7}
For $F$ and $T$ as in Theorem~\pr{flow1}, if $a\approx b$ then
$F_t(a) \approx F_t(b)$ for all $0 \leq t \leq T$.
\end{cor}

\begin{cor}\label{flow3}
For $F$ and $G$ as in Theorem~\pr{flow2}, if $F\equiv G$, then 
$F_t(a) \approx G_t(a)$ for all $0 \leq t \leq T$.
\end{cor}

\begin{pf}
By Proposition~\pr{18}(3) there is an infinitesimal $\beta$ for
the statement of Theorem~\pr{flow2}, which gives $\|F_t(a) - G_t(a) \| \leq  \beta t e^{Kt}$, so $F_t(a) \approx G_t(a)$.
\end{pf}

Our flow $F_t$ of a prevector field $F$ induces a classical flow on $M$ as follows.

\begin{dfn}\label{clflow}
Let $F$ be a  $D^1$ local prevector field. On a neighborhood $B'\su\rn$ and interval $[-T,T]$ as in Theorem~\pr{flow1}, we define the \emph{standard}
flow $h^F_t:B' \to M$ induced by $F$ as follows: $h^F_t(x) = st(F_t(x))$.
\end{dfn}

The following are immediate consequences of Theorems~\pr{flow1} and Corollary~\pr{flow3}.

\begin{thm}\label{hft}
Given a $D^1$ prevector field $F$ the following hold:
\begin{enumerate}
\item $h^F_t$ is Lipschitz continuous with constant $e^{K|t|}$.
\item $h^F_t$  is injective. 
\item If $G$ is another $D^1$ prevector field and $F \equiv G$ then $h^F_t = h^G_t$.
\end{enumerate}
\end{thm}

\begin{remark}\label{kf}
If $F$ is obtained from a classical vector field $X$ by the procedure of Example~\pr{e} then Keisler \ct[Theorem 14.1]{k} shows that our $h^F_t$ is in fact the flow of $X$ in the classical sense. By Theorem~\pr{hft}(3) this will be true for any prevector field $F$ that realizes $X$. 
\end{remark}

The results of this subsection have the following application to the standard setting.

\begin{clcor}\label{cc1}
For open $U \su\rn$ let $X,Y:U\to\rn$ be classical vector fields, where $X$ is Lipschitz with constant $K$, and $\| X(x)-Y(x) \| \leq b$ for all $x\in U$. If $x(t),x'(t)$ are integral curves of $X$ then 
$\|x(t)-x'(t)\| \leq  e^{Kt}\|x(0) - x'(0)\|$. If $y(t)$ is an integral curve of $Y$ with $x(0)=y(0)$ then 
$\|x(t) - y(t) \| \leq  {b\over K} (e^{Kt}-1)\leq b t e^{Kt}$.

\end{clcor}

\begin{pf}
Define prevector fields on $\s U$ by $F(a)-a=\DD X(a)$ and $G(a)-a=\DD Y(a)$ as in Example~\pr{e}, and apply Theorems \pr{flow1}, \pr{flow2}, and Remark~\pr{kf}.
\end{pf}

To conclude this section we look at the flow of a prevector field in an infinitesimal neighborhood of a fixed point. This corresponds to a zero of a vector field in the classical setting. In a neighborhood of such zero, one often approximates the given vector field with a simpler one (e.g. the linear approximation), to obtain an approximation of the original vector field's flow. We present the following approach for prevector fields,
which we apply in Section~\pr{secpend} to infinitesimal oscillations of a pendulum.

\begin{cor}\label{rsc}
Let $F,G$ be local $D^1$ prevector fields on $\s U$ where $U$ is a neighborhood of $p \in \RR^n$, and assume $F(p)=G(p)=p$. Fix
an infinitesimal $a>0$ and let $Q=\{x\in\s U : x-p\prec a\}$ (an external set).
\begin{enumerate}
 \item If      $F(x)-G(x) \prec\prec \DD a$ for all $x\in Q$ then $F_t(x)-G_t(x)\prec\prec a$ for all $x\in Q$, and the appropriate range of $t$, by which we mean finite $t$ for which $F_s(x) \in Q$ for all $0\leq s \leq t$.
\item If  $F(x)-G(x)  \prec\prec \|F(x)-x\|$ for all $x\in Q$,
then $F_t(x)-G_t(x) \prec\prec \|F_t(x)-p\|$ for all $ x\in Q$, and an appropriate range of $t$ as above.
\end{enumerate}
\end{cor}
\begin{pf}
For convenience assume $p=0$ so we have $F(0)=G(0)=0$ and $x\in Q$ simply means $x\prec a$. 
We first note that $F(x)-x =F(x)-x-F(0)+0 \prec \DD\|x\|$ by $D^1$ and Proposition~\pr{20}. Now
let $F'(x)=\frac{1}{a}F(ax)$, then $F'$ is defined for all $x\prec 1$ (i.e. $x\in \tth\RR^n$). 
For $x\prec 1$, $F'(x)-x=\frac{1}{a}\big(F(ax)-ax\big)
\prec\frac{1}{a}\DD\|ax\|=\DD\|x\|\prec\DD$ so $F'$ is a prevector field. We have $F'(x)-x-F'(y)+y=
\frac{1}{a}\big(F(ax)-ax-F(ay)+ay\big)\prec \frac{1}{a}\DD\|ax-ay\|=\DD\|x-y\|$ so $F'$ is $D^1$. 
Similarly define $G'$.

For (1) we have
$F'(x)-G'(x)=\frac{1}{a}\big(F(ax)-G(ax)\big)\prec\prec\frac{1}{a}\DD a=\DD$ so $F'\equiv G'$.
We thus get by Corollary~\pr{flow3} that $F'_t(x)\approx G'_t(x)$ for $x\prec 1$ and for appropriate range of $t$. Now, for $x\prec a$ we have 
$\frac{1}{a}x\prec 1$ so $F_t(x)-G_t(x)=a\big(F'_t(\frac{1}{a}x)-G'_t(\frac{1}{a}x)\big)\prec\prec a$.
(We remark that though our range of $t$ gives  $F'_s(x) \prec 1$ for all $0\leq s \leq t$, which is an external condition, in fact there is a finite ball $B$ such that $F'_s(x) \in \s B$ for all $0\leq s \leq t$.
This can be seen e.g. by underspill as in the proof of  Proposition~\pr{18}.)

For (2), the statement holds for $x=0$ since $0\prec\prec 0$.  Given a fixed $0 \neq x\prec a$ let $b=\|x\|$. Then for all $y\prec b$ we have 
$F(y)-G(y) \prec\prec \| F(y)-y\|  \prec\DD\|y\| \prec \DD b$, so  by (1) applied to $b$  we have $F_t(x)-G_t(x)\prec\prec b$ for appropriate range of $t$. By Theorem~\pr{flow1} we have $b= \|x-0\| \prec \|F_t(x)-F_t(0)\|=\|F_t(x)\|$, so together
$F_t(x)-G_t(x)\prec\prec  \|F_t(x)\|$.
\end{pf}

\begin{remark}\label{sc}
We can slightly weaken the assumptions in Corollary~\pr{rsc} by replacing the assumption that $G$ is $D^1$ 
by the weaker
assumption that all $x\prec a$ satisfy $G(x)-x \prec \DD\|x\|$. The proof remains unchanged. 
\end{remark}

\begin{example}\label{sc2}
In Corollary~\pr{rsc} we take $p=0\in\mathbb C$. We further assume that
0 is the only fixed point of $F,G$ in $Q$ (this corresponds to an isolated zero in the classical setting).
For $0\neq a,b\in\s\mathbb C$ we will say that $a$ and $b$ are \emph{adequal} if  $\frac{a}{b}\approx 1$.
Then Corollary~\pr{rsc}(2) tells us that if $F(x)-x$ and $G(x)-x$ are adequal for all $0\neq x\in Q$ then 
$F_t(x)$ and $G_t(x)$ are adequal for all $0\neq x\in Q$, and an appropriate range of $t$ as in Corollary~\pr{rsc}.  

\end{example}

\subsection{The canonical representative prevector field}

Once we have the standard function $h^F_t$, we can extend it to the nonstandard domain as usual, and 
use it to define a new prevector field $\tf$ as follows.
\begin{dfn}\label{can}
$\tf=h^F_\DD$. 
\end{dfn}

The map $\tf$ is indeed a prevector field, i.e. $\tf(a)-a \prec\DD$ for all $a$. Indeed, 
for $C \in\RR$ given by Proposition~\pr{18}(1) we have 
$\|F^n(a) - a\| \leq \sum_{m=1}^n \|F^m(a) - F^{m-1}(a)\| \leq nC\DD $, which implies
$\| h^F_t(a) -a\| \leq Ct$, which by transfer implies $\|\tf(a)-a\| \leq C\DD$.

 By Theorem~\pr{hft}(3), if  $F \equiv G$ then $\tf = \tg$. We will show in Theorem~\pr{tld} that $\tf \equiv F$, and so $\tf$ is a canonical choice
of a representative from the equivalence class of $F$.  (Perhaps in a smaller neighborhood of a given point, as required by Theorem~\pr{flow1}.) 
We will show in Propositions \pr{111}, \pr{112} that if $F$ is $D^1$ (resp. $D^2$) then $\tf$ is $D^1$ (resp. $D^2$).
That is, if a given equivalence class contains some member which is $D^1$ (resp. $D^2$) then the canonical representative $\tf$ of that equivalence class  is also $D^1$ (resp. $D^2$). 
We note that indeed not all members of the given class are $D^1$ (resp. $D^2$),
for example for $\s\RR$ take $F(x)=x$ for all $x\in\s\RR$, and $G(x)=x$ for all $x\neq 0$ and $G(0)=\DD^2$.
Then $F$ is $D^2$, $F\equiv G$, but $G$ is not even $D^1$, as is seen by taking $a=0,b=\DD^2$.

\begin{lemma}\label{ek}
Let $F$ be a local $D^1$ prevector field defined on $\s U$. Assume  $$\|F(a)-F(b)-a+b\| \leq K\DD\|a-b\|$$ 
for all  $a,b \in \s U$.
Then the flow of $F$ satisfies:
$$\| F^n(a)-F^n(b)-a+b \| \leq K \DD n e^{K\DD n} \|a-b\|.$$
\end{lemma}

\begin{pf} We have
\begin{align*} \| F^n(a)-F^n(b)-a+b \| &\leq
\sum_{i=1}^n \| F^i(a) - F^i(b) -F^{i-1}(a) +F^{i-1}(b) \|  \\
&\leq \sum_{i=1}^n  K\DD \|F^{i-1}(a) - F^{i-1}(b) \|  \\&
\leq \sum_{i=1}^n  K\DD  (1+K\DD)^{i-1}  \|a - b \|  
\leq  K \DD n e^{K\DD n} \|a-b\|. \end{align*} 
The third inequality is by internal induction as in the proof of Theorem~\pr{flow1}.
\end{pf}

\begin{prop}\label{111}
If $F$ is $D^1$ then $\tf$ is $D^1$.
\end{prop}

\begin{pf}
Assume $\|F(a)-F(b)-a+b \| \leq K\DD\|a-b\|$ for all $a,b$ in some $\s\! B$ as in Proposition~\pr{20}. Then for $n = \lf t / \DD  \rf$ we have by Lemma~\pr{ek}
$ \| F^n(a)-F^n(b)-a+b \| \leq K \DD n e^{K\DD n} \|a-b\|  \leq Kte^{Kt}\|a-b\|.$
So for standard $a,b$ we have $\|h^F_t(a)-h^F_t(b)-a+b\| \leq Kte^{Kt}\|a-b\|$.
Extending back to the nonstandard domain and evaluating at $t=\DD$ we get, by transfer, $\| \tf(a) - \tf(b) -a +b \| \leq K\DD e^{K\DD}\|a-b\|$.
\end{pf}

\begin{prop}\label{112}
If $F$ is $D^2$ then $\tf$ is $D^2$.
\end{prop}

\begin{pf} 
Assume $\|\xt F(a)\| \leq K\DD\|v\| \|w\|$ for all $a,v,w$ in some $\s\! B$ as in Proposition~\pr{202}.
We prove by internal induction that $$\| \xt F^n (a) \|=\| F^n(a)-F^n(a+v)-F^n(a+w)+F^n(a+v+w)\| \leq
K\DD\sum_{i=n-1}^{2n-2} (1+K\DD)^i \|v\| \|w\|.$$ 
 Let $p=F^n(a)$, $x=F^n(a+v)-F^n(a)$, $y=F^n(a+w)-F^n(a)$.
Then $\|x\| \leq (1+K\DD)^n \|v\|$, $\|y\| \leq (1+K\DD)^n \|w\|$, and 
\begin{align*} \|F(p+x+y) - & F^{n+1}(a+v+w))\| \\&\leq  
(1+K\DD)\|p+x+y - F^n(a+v+w)\|  \\& = (1+K\DD) \|-F^n(a)+F^n(a+v)+F^n(a+w)-F^n(a+v+w)\|  \\&
\leq(1+K\DD)K\DD\sum_{i=n-1}^{2n-2} (1+K\DD)^i \|v\| \|w\|  = K\DD\sum_{i=n}^{2n-1} (1+K\DD)^i \|v\| \|w\|, \end{align*} by the induction hypothesis.
Now \begin{align*} \|& F^{n+1}(a)-F^{n+1}(a+v)-F^{n+1}(a+w)+F^{n+1}(a+v+w)\| \\ \leq&
\|F(p) - F(p+x) - F(p+y) + F(p+x+y)\| +\|F(p+x+y)-F^{n+1}(a+v+w)\|  \\
\leq& K\DD\|x\|\|y\|  + K\DD\sum_{i=n}^{2n-1} (1+K\DD)^i \|v\| \|w\|   \\
\leq& K\DD(1+K\DD)^{2n} \|v\| \|w\| + K\DD\sum_{i=n}^{2n-1} (1+K\DD)^i \|v\| \|w\|
=K\DD\sum_{i=n}^{2n} (1+K\DD)^i \|v\| \|w\|, \end{align*} which completes the induction.

So for $n=\lf t/\DD \rf$  we have \begin{multline*} \| F^n(a)-F^n(a+v)-F^n(a+w)+F^n(a+v+w)\|  \\ \leq K\DD\sum_{i=n-1}^{2n-2} (1+K\DD)^i \|v\| \|w\|  
\leq K\DD n e^{2Kt} \|v\| \|w\| \leq Kt e^{2Kt} \|v \| \|w\|. \end{multline*}  Thus for standard $a,v,w$ we have 
$$\| h^F_t(a)-h^F_t(a+v)-h^F_t(a+w)+h^F_t(a+v+w)\| \leq  Kt e^{2Kt} \|v \| \|w\|.$$
Extending back to the nonstandard domain and evaluating at $t=\DD$ we get: 
$$\| \tf(a)-\tf(a+v)-\tf(a+w)+\tf(a+v+w)\| \leq  K\DD e^{2K\DD} \|v \| \|w\|.$$
\end{pf}

Next we would like to prove that $\tf \equiv F$.
We first need two lemmas. 

\begin{lemma}\label{101}
Let $F$ be a local prevector field defined on $\s U$. Assume  
$\|F(a)-a\| \leq C\DD$ and
$\|F(a)-a-F(b)+b\| \leq K\DD\|a-b\|$   for all  $a,b \in \s U$.
Then the flow of $F$ satisfies: $\|F^n(a)-a - n(F(a)-a)\| \leq KC n^2 \DD^2$.
\end{lemma}

\begin{pf} 
We have 
 \begin{align*} \|F^n(a)-a &- n(F(a)-a)\| = \|\sum_{i=1}^n \Bigl(F^i(a)-F^{i-1}(a) - (F(a)-a)\Bigr)\| \\ &\leq 
 \sum_{i=1}^n \| F(F^{i-1}(a)) - F^{i-1}(a) - F(a)+a \| 
\leq \sum_{i=1}^n K\DD  \|F^{i-1}(a) - a \|  \\& \leq  \sum_{1 \leq j < i \leq n} K\DD  \|F^j(a) - F^{j-1}(a) \|
\leq n^2 K\DD  C\DD. \end{align*}
\end{pf}

\begin{lemma}\label{102}
Let $v\in\s \RR^n$ with  $v \prec \DD$,  let $g:[0,\infty)\to\RR^n$ be the standard function $g(t) = st(\lf t / \DD \rf v )$, and let 
$V = st(v/\DD)$. Then $g(t)=tV$ and the extension of $g$ back to the nonstandard domain satisfies
$g( \DD)  \equiv v$.
\end{lemma}

\begin{pf} Let $n = \lf t / \DD \rf $. Then
$$\|tV -nv\| \leq \|tV -t(v/\DD)\| + \|t(v/\DD) - n\DD(v/\DD) \| = t\|V-(v/\DD)\|  +   | t - n\DD |  \|v/\DD\| \prec\prec 1.$$ 
This shows that $g(t)=tV$. So $g(\DD) = \DD V$, and  we have  $\|\DD V -v\| = \DD \|V- v/\DD\| \prec\prec \DD$. 
\end{pf}

We are now ready to prove the following.

\begin{thm}\label{tld}
If $F$ is a local $D^1$ prevector field then $\tf \equiv F$.
\end{thm}

\begin{pf}
By Proposition \pr{111} and Lemma~\pr{103} it is enough to show that $\tf(a) \equiv F(a)$  for all \emph{standard} $a$. So, for standard $a$ 
let $g_t(a)=st\Bigl(a+\lf t / \DD \rf (F(a)-a) \Bigr)$. Letting $n=\lf t / \DD \rf$
 we have 
$$\| h^F_t(a)- g_t(a) \| 
=\|st(F^n(a)) - st\Bigl(a+n (F(a)-a) \Bigr)\| = st\|F^n(a) - a-n (F(a)-a)\|
\leq At^2$$ 
for some $A \in \RR$, by Lemma~\pr{101}.
Extending and evaluating at $t=\DD$ gives $\|\tf(a)-g_\DD(a)\| \leq A\DD^2 \prec\prec\DD$, i.e. $\tf(a) \equiv g_\DD(a)$.
Now $g_t(a)-a= st\Bigl(\lf t / \DD \rf (F(a)-a) \Bigr)$
so by Lemma~\pr{102} we have 
$g_\DD(a)-a \equiv F(a)-a$, so $g_\DD(a) \equiv F(a)$, and together we get $\tf(a)\equiv F(a)$.
\end{pf}

To conclude, $\tf$ is a canonically chosen representative from the equivalence class of $F$ (perhaps in a smaller neighborhood of a given point), and if $F$ is $D^1$ (resp. $D^2$) then $\tf$ is $D^1$ (resp. $D^2$).

\subsection{Infinitesimal oscillations of a pendulum}\label{secpend}

We now demonstrate and discuss some of the concepts and results of this section in relation to a concrete physical problem, that of small oscillations of a pendulum. (Compare Stroyan \ct{str}.)

Let $x$ denote the angle between a pendulum and the downward vertical direction. 
By considering the projection of the force of gravitation in the direction of motion, one obtains the equation of motion  $$m\ell\ddot{x} = -mg\sin x$$ where $m$ is the mass of the bob of the pendulum, $\ell$ is the length of its massless rod, and $g$ is the constant of gravity. Letting $\omega=\sqrt{g/\ell}$ we have    $\ddot{x} = -\omega^2\sin x$.
The initial condition of releasing the pendulum at angle $a$ is described by $x(0)=a$, $\dot{x}(0)=0$.
We replace this single second order differential equation with the system of two first order equations $\dot{x}=\omega y$, $\dot{y}=-\omega\sin x$, and initial condition $(x,y)=(a,0)$. 
The classical vector field corresponding to this system is $X(x,y)=(\omega y,-\omega\sin x)$.

We are interested in ``small'' oscillations in the classical setting, i.e. the limiting behavior when the parameter $a$ above tends to 0, and correspondingly, infinitesimal oscillations in the hyperreal setting, i.e. when $a$ is  infinitesimal. To this end, if $p_a(t)$ is the classical motion with initial angle $a$, we look at the motion rescaled
by the factor $a$, i.e. we look at  $\frac{p_a(t)}{a}$. This is the $x$ component of the flow of the rescaled vector field
$Y(x,y)=\frac{1}{a}X(ax,ay)=(\omega  y, -\omega \frac{\sin ax}{a})$. The initial condition $(a,0)$ for $X$ corresponds to initial condition $(1,0)$ for $Y$. We can incorporate the parameter $a$ into our manifold and look at the vector field $Z$ on $\RR^3$ given by $Z(x,y,a)=(\omega  y, -\omega \frac{\sin ax}{a},0)$,
and initial conidtion $(1,0,a)$. Note that $Z$ is well defined and analytic also for $a=0$ (indeed 
$\frac{\sin ax}{a} = x - \frac{a^2x^3}{3!} + \frac{a^4x^5}{5!} - \cdots$), and its value for $a=0$ is 
$Z(x,y,0)=(\omega y , -\omega x, 0)$. 
The classical flow for $a=0$ i.e. initial condition $(1,0,0)$ is $(\cos\omega t,\sin\omega t,0)$, and so by Classical Corollary \pr{cc1} we have  $\frac{p_a(t)}{a} \to \cos\omega t$ (in fact, uniformly on finite intervals).
It follows that for infinitesimal $a$, $\frac{p_a(t)}{a} \approx \cos\omega t$, for all finite $t$. 

The above computation was for the classical flow of a classical vector field, and was then extended to the nonstandard domain. But we may also view the flow itself as occurring in the nonstandard domain $\s\RR^2$, via the prevector field 
$F(x,y)=(x+\DD\omega  y, y-\DD\omega \sin x)$ with initial condition $(a,0)$. This is 
the prevector field obtained from our classical vector field $X$ by the procedure of Example~\pr{e}. 
After rescaling  as before, we have
the prevector field $G(x,y)=(x+\DD\omega  y, y-\DD\omega  \frac{\sin ax}{a})$ and initial condition $(1,0)$.
Define the prevector field $E(x,y)=(x+\DD\omega  y, y-\DD\omega x)$, then for infinitesimal $a$ we have 
$G\equiv E$ since $E(x,y)-G(x,y)=(0,\DD\omega  (\frac{\sin ax}{ax}-1)x)$ and $\frac{\sin ax}{ax}-1\prec\prec 1$.
Let us define another prevector field $H(x,y)=(x\cos\DD\omega +y\sin\DD\omega , -x\sin\DD\omega + y\cos\DD\omega )$,
then $H$ is clockwise rotation of the $xy$ plane by angle $\DD\omega$, so
$H_t(1,0) \approx (\cos\omega t, -\sin\omega t)$.  
We have  $\cos\DD\omega-1\prec\prec\DD\omega$    and   $\sin\DD\omega-\DD\omega\prec\prec\DD\omega$, so 
$E\equiv H$. We have $G\equiv E\equiv H$,
 so by Corollary~\pr{flow3}, since $E$ is evidently $D^1$, $G_t(1,0) \approx (\cos\omega t, -\sin\omega t)$.  (We have used arguments from the proof of Corollary~\pr{rsc} rather than quoting it.)
So finally, the $x$ component of $G_t(1,0)$ is $\approx \cos\omega t$ for any infinitesimal $a$, which means that the $x$ component of $\frac{F_t(a,0)}{a}$ is $\approx \cos\omega t$ for any infinitesimal $a$.
We may thus say the following.
\begin{cor}\label{pendcor}
 The motion of a pendulum with infinitesimal amplitude $a$ is practically harmonic motion, in the sense that if rescaled to appreciable size, it is infinitely close to standard harmonic motion, for all finite time. 
\end{cor}

Equivalently, one could say that the motion itself is harmonic with
the given infinitesimal amplitude $a$, with error which is infinitely
smaller than $a$.

\section{Realizing classical vector fields}\label{rl}

Given a classical vector field on a smooth manifold $M$, we seek a prevector field realizing it.
Using Example~\pr{e} we can do this only locally, while by Proposition~\pr{re} these local prevector  fields are compatible up to equivalence. This leads to the following definition.

\begin{dfn}\label{coh1}
A $D^1$ (resp. $D^2$) \emph{coherent family} of local prevector fields on $M$ is a family $\{ (F_\alpha , U_\alpha) \}_{\alpha \in J}$ where $\{ U_\alpha \}_{\alpha \in J}$ is an open covering of $M$, and each $F_\alpha$ is a local $D^1$ (resp. $D^2$) prevector field on $\s U_\alpha$, such that 
for $\alpha,\beta\in J$,  $F_\alpha|_{\s\! U_\alpha  \cap  \s\! U_\beta} \equiv F_\beta|_{\s\! U_\alpha  \cap  \s\! U_\beta}$.
\end{dfn}

\begin{dfn}\label{coh2}
A coherent family $\{ (G_\alpha , V_\alpha) \}_{\alpha \in K}$ is said to be a \emph{refinement} of $\{ (F_\alpha , U_\alpha) \}_{\alpha \in J}$, if for each
$\alpha \in K$ there is $\beta \in J$ such that $V_\alpha \su U_\beta$ and $G_\alpha \equiv F_\beta|_{\s\! V_\alpha}$. 
\end{dfn}

\begin{dfn}\label{coh3}
A refinement $\{ (G_\alpha , V_\alpha) \}_{\alpha \in K}$ of $\{ (F_\alpha , U_\alpha) \}_{\alpha \in J}$
is said to be a \emph{flowing} refinement if there are $0 < T_\alpha \in \RR$ for each $\alpha \in K$
 such that the flow $h^{G_\alpha}_t$ is defined
on $V_\alpha$ for $0 \leq t \leq T_\alpha$.
\end{dfn}
By Theorem~\pr{flow1} any $D^1$ coherent family of prevector fields has a $D^1$ flowing refinement.
By Theorem~\pr{hft}(3), if $V_\alpha \cap V_\beta \neq \es$ then $h^{G_\alpha}_t=h^{G_\beta}_t$
on $ V_\alpha \cap  V_\beta$ for $0 \leq t \leq \min \{T_\alpha, T_\beta \}$.
If we can choose a single $0 < T \in \RR$ which is good for all $\alpha \in K$, then we will say that the original
family $\{(F_\alpha, U_\alpha) \}_{\alpha \in J}$ is \emph{complete}. In that case we have a global well defined flow $h_t : M \to M$ for $0 \leq t \leq T$, and by iteration, for all $0\leq t \in\RR$.
Extending $h_t$ back to $\sm$, let $G = h_\DD$, then $G$ is a global prevector field. By Proposition~\pr{111}, $G$ is $D^1$ since $\{F_\alpha \}$ is $D^1$, and by Proposition~\pr{112}, if $\{F_\alpha \}$ is $D^2$ then $G$ is $D^2$.
By Theorem~\pr{tld} we have $G|_{\s U_\alpha} \equiv F_\alpha$ for all $\alpha \in J$.
We will call $G$ the \emph{globalization} of the complete coherent family $\{(F_\alpha ,U_\alpha)\}_{\alpha \in J}$.
By Theorem~\pr{hft}(3) if two complete coherent families have a common refinement, then they define the same flow 
$h_t : M \to M$, and so they have the same globalization.

We note that if $\{(F_\alpha,U_\alpha) \}_{\alpha \in J}$ has a \emph{finite} flowing refinement, i.e. a flowing refinement $\{ (G_\alpha , V_\alpha) \}_{\alpha \in J}$
for which $J$ is finite, then $\{(F_\alpha,U_\alpha) \}_{\alpha \in J}$ is clearly complete.

\begin{dfn}\label{comps}
A coherent family $\{(F_\alpha, U_\alpha) \}_{\alpha \in J}$ has compact support, 
if there is a compact $C \su M$ such that $\{(F_\alpha, U_\alpha) \}_{\alpha \in J} \cup \{(I,M-C)\}$ is coherent, (recall $I(a)=a$ for all $a$).
\end{dfn}

Clearly a coherent $D^1$ family with compact support has a  finite flowing refinement, so the following holds.

\begin{prop}\label{cscp} 
A coherent $D^1$ family with compact support is complete. 
\end{prop}

Given a classical vector field $X$ on $M$ of class $C^1$ or $C^2$, we would like to realize it by a global prevector field on $\sm$ of class $D^1$ or $D^2$ respectively. In Proposition~\pr{6} we have shown that this can be done locally. 
We now state and prove our global realization result. 

In the following proof we use our assumption that our nonstandard extension satisfies countable saturation. This means that for any sequence $\{A_n\}_{n\in\N}$ of internal sets such that $A_n \neq\es$ and 
$A_{n+1} \su A_n$ for all $n$, one has $\bigcap_{n\in\N} A_n \neq\es$. 

\begin{thm}\label{gl}
Let $X$ be a classical $C^1$ (resp. $C^2$) vector field on $M$. Then there is a $D^1$ (resp. $D^2$) global prevector field $F$ on $\sm$ that realizes $X$, where the value of $F$ in $\tm$ is canonically prescribed. If $X$ has compact support (in the classical sense), then the value of $F$ throughout $\sm$ is canonically prescribed, with $F(a)=a$ for $a \in \sm - \tm$.
\end{thm}

\begin{pf}
Assume first that $X$ has compact support. 
There is a family $U_\alpha$ of coordinate neighborhoods for $M$, on each of which $X$ is realized by $F_\alpha$
as in Example~\pr{e}, and by Proposition~\pr{re} the family $\{ (F_\alpha,U_\alpha) \}$ is coherent. 
By Proposition~\pr{6}, the family $\{ (F_\alpha,U_\alpha) \}$ is $D^1$ (resp. $D^2$) if $X$ is $C^1$ (resp. $C^2$).
The vector field $X$ having compact support $C \su M$ in the classical sense implies 
that $\{ (F_\alpha,U_\alpha) \}$ has compact support in the sense of Definition~\pr{comps}.
Thus by Proposition~\pr{cscp} it is complete, and let $F$ be its globalization. 
We first notice that the flow $h_t:M \to M$ which defines $F$ satisfies $h_t(a)=a$ for all  $a \in M-C$ and 
so by transfer $F(a) = h_\DD(a) = a$ for all $a \in \sm - \s C \supseteq \sm - \tm$, proving the concluding statement regarding $X$ with compact support.
Furthermore,
by Propositions \pr{111}, \pr{112}, 
$F$ is $D^1$ (resp. $D^2$) if 
$\{ (F_\alpha,U_\alpha) \}$ is $D^1$ (resp. $D^2$), which, as mentioned, holds 
if $X$ is $C^1$ (resp. $C^2$).
By Proposition~\pr{re} and Theorem~\pr{hft}(3) $F$
is uniquely determined by  $X$.
This completes the compact support case. 

If $X$ does not have compact support, we proceed 
using countable saturation of our nonstandard extension.
Let $\{U_n\}_{n\in\N}$ be a sequence of open sets in $M$ with $\overline{U_n}$ compact, $\overline{U_n} \su U_{n+1}$, and $\bigcup U_n =M$.
Let $f_n:M \to [0,1]$ be a sequence of smooth functions with compact support, such that 
$f_n|_{U_{n+1}} = 1$. 
Now let $G_n$ be the realization of $f_n X$ given by  the compact support case.
Let $A_n = \{ F \in \s{\mbox{Map}} (M) : F|_{\s U_n} =  G_n |_{\s U_n} \}$, then $A_n$ is nonempty for each $n$, since $G_n \in A_n$. We further have $A_{n+1} \su A_n$ since 
$G_{n+1}|_{\s U_n} = G_n|_{\s U_n}$, which is true since $f_{n+1}$ and $f_n$ are 
both 1 on $U_{n+1} \supseteq \overline{U_n}$ and so the same flow determines  $G_{n+1}|_{\s U_n}$ and $G_n|_{\s U_n}$. 
So, by countable saturation $\bigcap A_n \neq\es$. An $F$ in this intersection satisfies  $F \in \s{\mbox{Map}} (M)$, 
i.e. it is internal. Since  $\bigcup \s U_n = \tm$, $F$ realizes $X$. 
The restriction $F|_\tm$ is uniquely determined by  $X$, since 
 $F|_{\s U_n} =  G_n |_{\s U_n}$  is uniquely determined by  $X$, again since  
$f_n$ is 1 on $U_{n+1} \supseteq \overline{U_n}$.
\end{pf}

In the following example we demonstrate the need for $\{ U_n \}$ and $\{ f_n \}$ in the proof 
of Theorem~\pr{gl}, and the fact that the values of $F$ on $\sm - \tm$ may depend on the choice of 
$\{ U_n \}$, $\{ f_n \}$. 

\begin{example}\label{e1}
Let $M=\ii$, and let $X$ be the classical vector field on $M$ given by $X(x)=-1$ for all $x \in \ii$.
On $\s\ii$ $X$ does not induce a prevector field via the procedure of Example~\pr{e} since for $\DD > x \in\s\ii$, 
$x-\DD \not\in \s\ii$. However we can take the coherent family $\{(F_r , (r,1))\}_{r>0}$ where $F_r$ is always defined by $F_r(a)=a-\DD$. The standard flow $h^{F_r}_t$ is defined for $0 \leq t \leq r$ and always given by $h^{F_r}_t(a)=a-t$.
But this family is not complete. There is no common $T>0$ for which the flow is defined on $[0,T]$, and so there is no global flow $h_t:\ii\to\ii$ in which one can substitute $t=\DD$. (Note that the global prevector field that may seem to exist by naively substituting $t=\DD$ ignoring the problem of common domain $[0,T]$, would be $a \mapsto a-\DD$, which, as noted, is not defined on $\s\ii$.)

So, following the proof of Theorem~\pr{gl}, let $\{a_n\}$ be a strictly decreasing sequence with $a_n \to 0$.
Let $U_n = (a_n , 1-a_n)$ and let $f_n:\ii \to [0,1]$ be a smooth function such that $f_n(x)=1$ 
for $a_{n+1}   \leq x \leq 1-a_{n+1}$, and $f_n(x)=0$ for $ 0<x \leq a_{n+2}$ 
and $1-a_{n+2} \leq x < 1$. 
To realize $f_n X$ as in Example~\pr{e} we do not need a covering $\{ U_\alpha \}$ as in the general case
appearing in the proof of Theorem~\pr{gl}, 
rather we can take one $F_n$ defined on all $\ii$.
For $a \in \s (a_{n+1}, 1-a_{n+1})$ we have  $F_n(a)=a-\DD$, and so 
for $a \in (a_n , 1-a_n)$ and $0 \leq t \leq a_n-a_{n+1}$
we have  $h^{F_n}_t(a)=a-t$, and so finally
for $a \in \s U_n = \s (a_n , 1-a_n)$ the realization $G_n$ of
$f_n X$ satisfies $G_n(a)=a - \DD$.
Thus, a global $F:\s\ii \to\s\ii$ which is obtained from the sequence $G_n$ as in the proof of Theorem~\pr{gl} will have $F(a)=a-\DD$ for all $a \in {}^\ha {\ii} = \{ a \in \s \ii \ : 0 < \ st(a) < 1 \}$, and this fact
is independent of all choices involved in the construction. 
However, the values on $\s\ii - {}^\ha {\ii}$ may indeed  depend on our choice of $\{U_n\}$ and $\{f_n\}$, as we now demonstrate.

Suppose our nonstandard extension is given by the ultrapower construction on the index set $\N$ with nonprincipal ultrafilter,
and elements in the ultrapower are   given by sequences in angle brackets $\langle x_i \rangle_{i \in\N}$.
\footnote{Such extension always satisfies countable saturation.}
Assume $\DD= \langle \delta_i \rangle_{i \in\N}$ where $\{\delta_i\}$ is  a strictly decreasing sequence 
with $\delta_i \to 0$. Then $G_n = h^{F_n}_\DD = \langle h^{F_n}_{\delta_i} \rangle_{i\in\N}$.
Let $F=\langle h^{F_i}_{\delta_i} \rangle_{i\in\N}$, and we claim that $F|_{\s U_n} = G_n|_{\s U_n}$
for all $n$, i.e. $F \in \bigcap A_n$. 
Indeed, the elements of $\s U_n$ are represented by sequences $\langle u_i \rangle_{i \in\N}$ 
such that $u_i \in U_n$ for all $i$, 
and so for $i$ sufficiently large so that $i \geq n$ and $\delta_i < a_n -a_{n+1}$ we have 
$h^{F_i}_{\delta_i}(u_i) = u_i -\delta_i = h^{F_n}_{\delta_i}(u_i)$.
Now let $x=\langle a_{i+2} \rangle_{i\in\N}$  and $y=\langle a_i + \delta_i \rangle_{i\in\N}$, then 
 $F(x)=x$ and $F(y)=y-\DD$. 
If we repeat our construction with $a'_n = a_{n-2}+\delta_{n-2}$ in place of $a_n$, producing the realization $F'$, then for the same reason that $F(x)=x$ we will have $F'(y)=y \neq F(y)$, showing that
$F$ indeed depends on our choices.  
\end{example}

\section{Lie bracket}\label{lbr}

Given two local prevector fields $F,G$ for which $F^{-1},G^{-1}$ exist, e.g. if $F,G$ are $D^1$ (by Remark~\pr{finv2}),
we define their Lie bracket $[F,G]$ 
 as follows. Its relation to the classical Lie bracket will be clarified in Section~\pr{rclb}.
\begin{dfn}\label{dflie}
$[F,G] = (G^{-1} \circ F^{-1} \circ G \circ F)^{\lf {1 \over \DD}   \rf}_{\phantom i}$. 
\end{dfn}
Since our  fixed choice of $\DD$ was arbitrary, we may have chosen it as $1 \over N$ for some infinite  $N \in \s\N$, and so we may assume $1\over \DD$ is in fact a hyperinteger and drop the $\lf \cdot \rf$ from the above expression.
In Theorem~\pr{rb} below we will justify this definition, i.e. we will establish its relation to the classical Lie bracket.
We will show that if $F,G$ are  $D^1$ then $[F,G]$ is indeed a prevector field,  and if $F,G$ are $D^2$ then $[F,G]$ is $D^1$.
Furthermore, we will show that if $F,G$ are  $D^2$ and $F \equiv F'$, $G \equiv G'$ then $[F,G] \equiv [F',G']$. 
We will give an example showing that this is not true if $F,G$ are merely  $D^1$. 
We will show that the Lie bracket of two $D^2$ prevector fields is equivalent to the identity prevector field if and only if their local standard flows commute. 
In the present section our study will always be local, and so the quantifier ``for all $a$'' will always mean for all $a$ in $\s U$ where $U$ is some appropriate coordinate neighborhood, and all computations are in coordinates.

\subsection{Fundamental properties of Lie bracket}

\begin{thm}\label{7}
If $F,G$ are  local $D^1$ prevector fields then $[F,G]$ is a prevector field, that is, $[F,G](a)-a \prec\DD$ for all $a$.
\end{thm}

\begin{pf} 
Substituting $x=a$ and $y=F^{-1} \circ G \circ F(a)$ 
in the relation $F(x)-x-F(y)+y \prec\DD\|x-y\|$
gives 
$$F(a)-a-G \circ F(a) + F^{-1} \circ G \circ F(a) \prec \DD \|a-F^{-1} \circ G \circ F(a)\| \prec \DD^2.$$
Now substituting $x=F(a)$ and $y=G^{-1} \circ F^{-1} \circ G \circ F(a)$ 
in the relation $G(x)-x-G(y)+y \prec\DD\|x-y\|$
 gives
 $$G\circ F(a)-F(a)- F^{-1} \circ G \circ F(a)+ G^{-1} \circ F^{-1} \circ G \circ F(a) \prec \DD \|F(a)-G^{-1} \circ F^{-1} \circ G \circ F(a)\| \prec \DD^2.$$
Adding the above two expressions gives: $ G^{-1} \circ F^{-1} \circ G \circ F(a)  - a \prec  \DD^2$.
By underspill in an appropriate $\s U$ there exists $C \prec 1$ such that  $\| G^{-1} \circ F^{-1} \circ G \circ F(a)  - a \|\leq C \DD^2$ for all $a\in\s U$. And so 
 \begin{multline*}\| (G^{-1} \circ F^{-1} \circ G \circ F)\1 (a)  - a \| \\ \leq \sum_{k=1}^{1\over \DD}
\| (G^{-1} \circ F^{-1} \circ G \circ F)^k (a)  - (G^{-1} \circ F^{-1} \circ G \circ F)^{k-1} (a) \|
\leq C \DD.\end{multline*}
\end{pf}

\begin{example}\label{e2}
We give an example of two prevector fields $F,G$, where $F$ is $D^2$ (so also $D^1$) and $[F,G]$ is not a prevector field.
Let $M=\RR^2$ and let $F(x,y)=(x+\DD,y)$, $G(x,y)= (x, y+\DD\sin{\pi \over 2\DD}x)$. 
Then   $[F,G](0,0) = (0,1)$ so $[F,G](0,0)-(0,0) = (0,1) \not\prec \DD$ .
\end{example}

To prove that if $F,G$ are  $D^2$ then $[F,G]$ is $D^1$ we  need the following  lemma. 
A sum of eight terms appears in
 its statement, namely 
$$\Bigl(F(a)-a\Bigr)-\Bigl(F(b)-b\Bigr)  -\Bigl(F(G(a))-G(a)\Bigr) +  \Bigl(F(G(b))-G(b)\Bigr)$$
which is similar to the sum
\begin{multline*}
\xt (F-I)(a)= \\
\Bigl(F(a)-a\Bigr)-\Bigl(F(a+v)-(a+v)\Bigr)  -\Bigl(F(a+w)-(a+w)\Bigr) +  \Bigl(F(a+v+w)-(a+v+w)\Bigr)
\end{multline*}
appearing in the general definition of $D^k$ applied to $k=2$.
As already noticed, the four terms $a,a+v,a+w,a+v+w$ cancel, leaving the four terms appearing in 
Definition~\pr{d2}.   
In the present sum the corresponding four terms
$a,b,G(a),G(b)$ do not cancel, and we remain with all eight terms.
We have already encountered a similar eight term sum  $\xt(\f\circ F -\f)(a)=  \xt(G \circ \f -\f)(a)$
where no cancellation occurs, in the proof of Proposition~\pr{4}.

\begin{lemma}\label{8}
Let $F$ be $D^2$  and  $G$ be $D^1$, then for all $a,b$ with $a-b \prec \DD$,
$$F(a)-F(b)-F(G(a))+F(G(b))-a+b+G(a)-G(b) \prec \DD^2\|a-b\|.$$
\end{lemma}

\begin{pf}
Let  $v=b-a$ and $w=G(a)-a$. 
Since  $F$ is  $D^1$  (by Proposition~\pr{3}) we have 
$F(a+v+w)-F(G(b))-(a+v+w)+G(b) \prec \DD\|a+v+w-G(b)\|$. But $a+v+w = b+G(a)-a$ and so we have 
$$F(a+v+w)-F(G(b))-b-G(a)+a+G(b) \prec \DD \| b+G(a)-a-G(b) \| \prec \DD^2 \|a-b\|$$
since $G$ is $D^1$. So 
\begin{multline*}
\|F(a)-F(b)-F(G(a))+F(G(b))-a+b+G(a)-G(b) \|   \\ = 
 \|F(a)-F(a+v)-F(a+w)+F(a+v+w)-F(a+v+w) + F(G(b))-a+b+G(a)-G(b) \| \\
\leq \|  F(a)-F(a+v)-F(a+w)+F(a+v+w)     \|  + \|      -F(a+v+w) + F(G(b))-a+b+G(a)-G(b)        \| \\
\prec \DD\|b-a\|\|G(a)-a\| + \DD^2\|a-b\| \prec \DD^2\|a-b\|.
\end{multline*}
\end{pf}

\begin{thm}\label{9}
If $F,G$ are $D^2$  then $[F,G]$ is $ D^1$.
\end{thm}

\begin{pf}
By Propositions \pr{3}, \pr{30}, and \pr{40},  $F^{-1} \circ G \circ F$ is  $D^1$.
Now in Lemma~\pr{8} take $G$ to be $F^{-1} \circ G \circ F$ then we  get for $a-b \prec \DD$:
$$F(a)-F(b)-G\circ F(a) + G\circ F(b)-a+b+F^{-1} \circ G \circ F (a) - F^{-1} \circ G \circ F (b)  \prec \DD^2 \|a-b\|.$$
As above $G^{-1} \circ F^{-1} \circ G$ is  $D^1$ 
and now take in Lemma~\pr{8} $a,b,F,G$ to be respectively $F(a),F(b),G,G^{-1} \circ F^{-1} \circ G$ then we get

\begin{multline*}
G \circ F(a)- G \circ F(b)-F^{-1} \circ G\circ F(a) + F^{-1} \circ G\circ F(b)
 \\ -F(a)+F(b)+G^{-1} \circ F^{-1} \circ G \circ F (a) - G^{-1} \circ F^{-1} \circ G \circ F (b)  \prec \DD^2 \|F(a)-F(b)\| \prec \DD^2\|a-b\|\end{multline*} by Proposition~\pr{0}.
Adding these two inequalities we get 
$$G^{-1} \circ F^{-1} \circ G \circ F (a) - G^{-1} \circ F^{-1} \circ G \circ F (b) -a+b  \prec \DD^2\|a-b\|$$

Denote $H=G^{-1} \circ F^{-1} \circ G \circ F$ then $[F,G] = H\1 $ and 
so we must show $H\1 (a)-H\1 (b)-a+b \prec \DD\|a-b\|$ and we know 
$H(a)-H(b)-a+b \prec \DD^2\|a-b\|$. By underspill in an appropriate $\s U$ there exists $C \prec 1$ such that
$\|H(a)-H(b)-a+b \| \leq C\DD^2\|a-b\|$ for all $a,b \in\s U$.
So by Lemma~\pr{ek} with $K=C\DD$ and $n={1\over \DD}$, we get
$\|H\1 (a)-H\1 (b)-a+b \| \leq C\DD e^{C\DD}\|a-b \|$.
\end{pf}

\begin{example}\label{e3}
We give an example of two prevector fields $F,G$, where $F$ is $D^2$, $G$ is  $D^1$ and $[F,G]$ is not $D^1$.
Let $M=\RR^2$ and let $F(x,y)=(x+\DD,y)$, $G(x,y)= (x, y+\DD^2 \sin{\pi \over 2\DD}x)$. 
Clearly $F$ is $D^2$, and we show $G$ is $D^1$: 
\begin{align*}
\|G(x_1,y_1) -(x_1,y_1) - G(x_2,y_2)+(x_2,y_2)\| 
&= \DD^2 |\sin{\pi \over 2\DD}x_1 - \sin{\pi \over 2\DD}x_2 | \\ &=  \DD{\pi \over 2}|(x_1 - x_2) \cos{\pi \over 2\DD}\theta |   \\&
\prec \DD|x_1 - x_2| \prec \DD\|(x_1,y_1)-(x_2,y_2)\|,\end{align*} 
where $x_1 \leq \theta \leq x_2$. 
Finally we show $[F,G]$ is not $D^1$:   $[F,G](0,0) = (0,\DD)$,  $[F,G](\DD,0) = (\DD,-\DD)$, so $[F,G](0,0)-(0,0) - [F,G](\DD,0)+(\DD,0)  
= (0,2\DD) \not\prec \DD\|(0,0)-(\DD,0)\|$.

\end{example}

Our definition of Lie bracket involves an iteration $1\over \DD$ times of the commutator
 $G^{-1} \circ F^{-1} \circ G \circ F$. The following Proposition compares this with multiplication by $1\over \DD$ in coordinates. It will be used in the proofs of Theorems \pr{ld2}, \pr{rb}, \pr{comm}.

\begin{prop}\label{scl}
Let $F,G$ be $D^2$,
then $[F,G] (a)\equiv a+{1 \over \DD} \Bigl(G^{-1} \circ F^{-1} \circ G \circ F(a) - a\Bigr)$  for all $a$.
\end{prop}

\begin{pf}
Let $H=G^{-1} \circ F^{-1} \circ G \circ F$.
 The proof of Theorems \pr{7} provides $C' \prec 1$ such that 
$\|H(a)-a\| \leq C' \DD^2$ for all $a$. The proof of Theorem~\pr{9} provides $C''\prec 1$
such that $\|H(a)-H(b)-a+b\|\leq C'' \DD^2\|a-b\|$ for all $a,b$.
Taking $C=C'\DD$, $K=C''\DD$ and $n={1\over\DD}$ in Lemma~\pr{101} we get 
$\|H\1 (a)-a - {1\over\DD}(H(a)-a)\| \leq C'C'' \DD^2\prec\prec\DD$.
\end{pf}

Next we would like to show that if $F,F',G,G'$ are  $D^2$ and $F \equiv F'$, $G \equiv G'$ then $[F,G] \equiv [F',G']$. We will need the following two lemmas.

\begin{lemma}\label{10}
If $G,H$ are  $D^2$ and $G \equiv H$, (i.e. $G(a)-H(a)\prec\prec\DD$ for all $a$) then 
$$\Bigl(G(a)-H(a)\Bigr)-\Bigl(G(b)-H(b)\Bigr)\prec\prec\DD\|a-b\|$$ 
for all $a,b$ with $a-b \prec \DD$.
\end{lemma}

\begin{pf}
Let $F(x)=G(x)-H(x)$ so $F(x) \prec\prec\DD$ for all $x$. Assume $F(a)-F(b)$ is not $\prec\prec \DD\|a-b\|$ 
for some $a,b$ with $a-b \prec \DD$,
then $ \DD\|a-b\| \prec\| F(a)-F(b)\|$.
Let $v=b-a$ then $\DD\|v\|^2 \prec \|v \| \|F(a)-F(a+v)\|$, and since
$G,H$ are $D^2$, $F$ satisfies
$\x^2_{v,v}F(x)  \prec \DD\|v\|^2$ for all $x$. 
Together we have 
  $\x^2_{v,v}F(x) \prec \|v\| \|F(a)-F(a+v)\|$,
so by Lemma~\pr{1} (taking some ball around $st(a)$),
there is $m \in\s\N$ such that 
$F(a)-F(a+v)\prec \|v\| \| F(a)-F(a+mv)\| \leq \|v\|  \Bigl(\| F(a)\|+\|F(a+mv)\|\Bigr) \prec\prec \|v\|\DD$.
\end{pf}

\begin{lemma}\label{11}
If $G,H$ are prevector fields with  $G\equiv H$ and $G$ is  $D^1$, then $G^{-1}\equiv H^{-1}$ (assuming $H^{-1}$ exists).
\end{lemma}

\begin{pf}
Given $a$ let $x=G^{-1}(a)$ and $y=H^{-1}(a)$ then we must show $x-y \prec\prec\DD$.
We have $G(x)=a=H(y)$ so $\|G(x)-G(y)\| = \|H(y)-G(y)\| =\beta\DD$ 
for some $\beta \prec\prec 1$. Since $G$ is  $D^1$, 
$\|G(x)-G(y) - x+y\| = K\DD\|x-y\|$ for some $K \prec 1$. So
$$\|x-y\| \leq \|G(x)-G(y) - x+y\|  + \|G(x)-G(y)\| =  K\DD\|x-y\| + \beta\DD.$$ 
So $(1-K\DD)\|x-y\| \leq \beta\DD$ or $\|x-y\| \leq {\beta \over 1-K\DD}\DD \prec\prec\DD$.
\end{pf}

We are now ready to prove the following.

\begin{thm}\label{ld2}
 If $F,E,G,H$ are  $D^2$,  $F \equiv E$ and $G \equiv H$ then $[F,G] \equiv [E,H]$.
\end{thm}

\begin{pf}
We first claim that it is enough to establish the statement with $F=E$, that is, to show $[F,G] \equiv [F,H]$. 
Indeed it is clear from the definition that $[G,F]=[F,G]^{-1}$, so if we know  $[F,G] \equiv [F,H]$ and similarly 
$[H,F] \equiv [H,E]$ then by Theorem~\pr{9} and Lemma~\pr{11}   we have
$[F,G] \equiv [F,H]=[H,F]^{-1} \equiv [H,E]^{-1}=[E,H] $. 

So we proceed to show  $[F,G] \equiv [F,H]$. Given $x$ let $a=G( F(x))$, $b=H( F(x))$, 
then by assumption $a-b \prec\prec\DD$. 
By Propositions \pr{3}, \pr{30}, $$F^{-1}(a)-F^{-1}(b)-a+b \prec \DD\|a-b\|\prec\prec \DD^2.$$ 
Denote $c=G^{-1} \circ F^{-1} \circ G \circ F(x)$.
By Lemma~\pr{10} 
$$F^{-1}(a)-H(c)-a+b = \Bigl(G(c) - H(c) \Bigr) - \Bigl(    G(F(x)) - H(F(x))    \Bigr)
\prec\prec \DD \| c-F(x)\| \prec \DD^2.$$
Combining the last two inequalities we get $H(c) - F^{-1}(b) \prec\prec\DD^2$ and so
\begin{multline*}G^{-1} \circ F^{-1} \circ G \circ F(x) - H^{-1} \circ F^{-1} \circ H \circ F(x)  
=\\  H^{-1}(H(c)) - H^{-1}(F^{-1}(b))   \prec \| H(c) - F^{-1}(b) \|  \prec\prec  \DD^2\end{multline*}
by Propositions \pr{3}, \pr{30}, \pr{0}.
So we have $${1\over \DD}\Bigl(G^{-1} \circ F^{-1} \circ G \circ F(x)\Bigr) - {1\over \DD}\Bigl(H^{-1} \circ F^{-1} \circ H \circ F(x)\Bigr) \prec\prec\DD$$
and so by Proposition~\pr{scl} $[F,G](x)\equiv [F,H](x)$.
\end{pf}

\begin{example}\label{e4}
We  give an example of $F,H$ which are  $D^2$, $G$ is $D^1$ and
 $G \equiv H$, and yet $[F,G] \not\equiv [F,H]$.
Let $M=\RR^2$ and let $F(x,y)=(x+\DD,y)$, $H(x,y)=(x,y)$, and $G(x,y) = (x, y+\DD^2\sin{\pi \over 2\DD}x)$. 
Then $[F,H](0,0)=(0,0)$ whereas  $[F,G](0,0) = (0,\DD) \not\equiv (0,0)$.
Clearly $F,H$ are  $D^2$, and it has been shown in Example~\pr{e3} that $G$ is $D^1$.
\end{example}

\subsection{Relation to classical Lie bracket}\label{rclb}

The following theorem justifies our definition of $[F,G]$, by relating it to the classical notion of Lie bracket.

\begin{thm}\label{rb}
Let $X,Y$ be two classical $C^2$ vector fields and let $[X,Y]_{cl}$ denote their classical Lie bracket. 
Let $F,G$ be $D^2$ prevector fields that realize $X,Y$ respectively.
Then $[F,G]$ realizes $[X,Y]_{cl}$. 
\end{thm}

\begin{pf}
By Remark~\pr{kf}, the flows $h^F_t$, $h^G_t$ coincide with the classical flows of $X$, $Y$. 
It is well known that $[X,Y]_{cl}$ is related in coordinates to the classical flow as follows:
$$[X,Y]_{cl}(p)=\lim_{t\to 0} {1 \over t^2} \Bigl((h^G_t)^{-1} \circ  (g^F_t)^{-1} \circ h^G_t \circ h^F_t(p) - p\Bigr).$$ 
By the equivalent characterization of limits via infinitesimals  we thus have  
$$[X,Y]_{cl}(p) \approx {1 \over \DD^2} \Bigl(\tg^{-1} \circ  \tf^{-1} \circ \tg \circ \tf(p) - p\Bigr).$$ 
Now, if $v \approx w$ then $\DD v \equiv \DD w$, so by Example~\pr{e}, $[X,Y]_{cl}$ can be realized by the prevector field 
$$A(a) =  a+{1 \over \DD} \Bigl(\tg^{-1} \circ  \tf^{-1} \circ \tg \circ \tf(a) - a\Bigr).$$ Thus it remains to show that
$[F,G]  \equiv A$. 
By Proposition~\pr{112} $\tf,\tg$ are $D^2$, and so 
by Proposition~\pr{scl} $ [\tf,\tg]  \equiv A $.
By Theorem~\pr{tld} $F\equiv \tf, G\equiv \tg$, 
and so by Theorem~\pr{ld2} $ [F,G]  \equiv A $.
\end{pf}

The following theorem corresponds to the classical fact that the bracket of two vector fields vanishes if and only if their flows commute.

\begin{thm}\label{comm}
Let $F,G$ be two $D^2$ prevector fields. Then $[F,G]\equiv I$ (recall $I(a)=a$ for all $a$), if and only if 
$h^F_t \circ h^G_s  =  h^G_s \circ h^F_t $ for all $0 \leq t,s \leq T$ for some $0 < T \in\RR$. 
\end{thm}

\begin{pf}
Assume first that $[F,G]\equiv I$, i.e. $[F,G](a)- a\prec\prec\DD$ for all $a$.  
So by Proposition~\pr{scl} ${1 \over \DD} \Bigl(G^{-1} \circ F^{-1} \circ G \circ F(a) - a\Bigr) \prec\prec\DD$,
so  $G^{-1} \circ F^{-1} \circ G \circ F(a) - a \prec\prec\DD^2$, which implies by Proposition~\pr{0} that
$G \circ F(a) - F\circ G( a) \prec\prec\DD^2$ for all $a$.
Now let $n = \lf t/\DD \rf$ and $m = \lf s/\DD \rf$, then we need to show 
$F^n \circ G^m(a) \approx G^m \circ F^n(a)$ for all $a$.
This involves $nm$ interchanges of $F$ and $G$, where a typical move is from
$F^k \circ G^r \circ F \circ G^{m-r} \circ F^{n-k-1}$ to $F^k \circ G^{r+1} \circ F \circ G^{m-r-1} \circ F^{n-k-1}$. Applying $F \circ G(p) - G\circ F( p) \prec\prec\DD^2$ to $p=G^{m-r-1} \circ F^{n-k-1}(a)$
we get $$F \circ G^{m-r} \circ F^{n-k-1}(a) -G \circ F \circ G^{m-r-1} \circ F^{n-k-1}(a) \prec\prec\DD^2.$$ 
By Propositions \pr{3}, \pr{20} there is $K\in\RR$ such that
$\| F(a)-a-F(b)+b\| \leq K\DD\|a-b\|$  and $\| G(a)-a-G(b)+b\| \leq K\DD\|a-b\|$ for all $a,b$ in an appropriate domain.
Then by Theorem~\pr{flow1} applied to $G^r$ and then to $F^k$,
\begin{multline*}
\|F^k \circ G^r \circ F \circ G^{m-r} \circ F^{n-k-1}(a)-F^k \circ G^{r+1} \circ F \circ G^{m-r-1} \circ F^{n-k-1}(a) \|  \\
\leq e^{K(t+s)}\| F \circ G^{m-r} \circ F^{n-k-1}(a) -G \circ F \circ G^{m-r-1} \circ F^{n-k-1}(a)\| \prec\prec\DD^2.
\end{multline*}
Adding the $nm$ contributions when passing from $F^n \circ G^m(a)$ to $G^m \circ F^n(a)$
we get $$F^n \circ G^m(a) - G^m \circ F^n(a) \prec\prec 1.$$
This is because among the $nm$ differences that we add, there is a maximal one, which is say $\beta \DD^2$ with $\beta \prec\prec 1$, and so the sum of all $nm$ contributions is $\leq nm \beta\DD^2 \leq ts\beta \prec\prec 1$.

Conversely, assume $h^F_t \circ h^G_t  =  h^G_t \circ h^F_t $. 
Then by transfer $\tf \circ \tg = \tg \circ \tf$, so $\tg^{-1} \circ  \tf^{-1} \circ \tg \circ \tf = I$, and so 
$ [\tf,\tg]=I$. By Proposition~\pr{112} and Theorems \pr{tld}, \pr{ld2} we get $[F,G]\equiv I$.
\end{pf}

We have the following application to the standard setting.

\begin{clcor}\label{cc4}
Let $X,Y$ be classical $C^2$ vector fields. Then the flows of $X$ and $Y$ commute if and only if their Lie bracket vanishes. 

It follows that if $X_1,\dots,X_k$ are $k$ independent vector fields with $[X_i,X_j]_{cl}=0$ (classical Lie bracket) for $1 \leq i,j \leq k$, then there are coordinates in a neighborhood of any given point such that
$X_1,\dots,X_k$ are the first $k$ coordinate vector fields.
\end{clcor}

\begin{pf}
Define prevector fields by $F(a)=a+\DD X(a)$ and $G(a)=a+\DD Y(a)$ as in Example~\pr{e}, and apply 
Proposition~\pr{6}, Remark~\pr{kf}, and  Theorems \pr{rb}, \pr{comm}. The final statement is a 
straightforward conclusion in the classical setting.
\end{pf}

\end{document}